\documentclass[a4paper]{article}

\usepackage{graphicx}
\usepackage{float}
\usepackage{subfigure}
\usepackage{amsfonts}
\usepackage{amsmath}
\usepackage{amsthm}
\usepackage{amssymb}
\usepackage{bm}
\usepackage{color}
\usepackage[T1]{fontenc}
\usepackage[utf8]{inputenc}
\usepackage{comment}
\usepackage{authblk}
\usepackage{multirow}
\usepackage{mathtools}
\usepackage{lineno}
\usepackage[maxbibnames=9]{biblatex} 
\newenvironment{system}%
	{\left\lbrace \begin{array}{@{} l @{} }}%
	{ \end{array}\right.}

\newtheorem{lem}{Lemma}
\newtheorem{thm}{Theorem}

\newtheorem{remark}{Remark}

\newcommand{\bbR}{\mathbb{R}}
\newcommand{\bbN}{\mathbb{N}}

\newcommand{\calO}{\mathcal{O}}
\newcommand{\calV}{\mathcal{V}}

\newcommand{\calX}{\mathcal{X}}

\DeclareMathOperator{\ctg}{ctg}
\DeclareMathOperator{\cst}{cst}

\title{A barycentric trigonometric Hermite interpolant via an iterative approach} 

\author[a]{Giacomo Elefante}
\affil[a]{Department of Mathematics \lq\lq Tullio Levi-Civita\rq\rq, University of Padova, Via Trieste 63, Padova, Italy\\ \texttt{giacomo.elefante@unipd.it}}

\date{}

\begin{document}

\maketitle

\begin{abstract} 	
In this work we construct an Hermite interpolant starting from basis functions that satisfy a Lagrange property. In fact, we extend and generalise an iterative approach, introduced by Cirillo and Hormann \cite{CH18} for the Floater-Hormann family of interpolants. Secondly, we apply this scheme to produce an effective barycentric rational trigonometric Hermite interpolant at general ordered nodes using as basis functions the ones of the trigonometric interpolant introduced by Berrut \cite{Berrut88}. For an easy computational construction, we calculate analytically the differentation matrix. Finally, we conclude with various examples and a numerical study of the convergence at equidistant nodes and conformally mapped nodes.
\end{abstract}

\section{Introduction}

Hermite interpolation is a classical problem which has been widely explored. It has many applications, and therefore, it is still vivid and studied. 

The problem consists into finding a good function that satisfies the Hermite conditions, i.e.,
\begin{equation} \label{HermCond}
    r_m^{(j)}(x_i)=f_{i,j}=f^{(j)}(x_i),\quad i=0,\dots,n;\,\,\, j=0,\dots,m ,
\end{equation}
which means that it has to interpolate not only a function on some nodes $x_i, i=0,\dots,n$ but also the values of its first $m$ derivatives.  

A good proposal for a solution of this problem could be given by interpolants based on rational functions since they may take advantage of the well-known stability of the latters for the classical problem of interpolation (see e.g. \cite{Schneider91}). In fact, some recent works go in that direction, as the family of Hermite barycentric interpolants introduced in \cite{Jing20}, the family of interpolants in \cite{CH18} or the combination of rational functions and six-point Shepard basis functions introduced in \cite{DellAccio20}.


A good stable interpolant that achieves a fast and accurate solution is preferable since it may be useful for the many applications where an Hermite interpolant is used, such as in Treecode algorithm \cite{Krasny19} or constructing curves and surfaces for images \cite{Oumellal21} or approximating the reliability of a hammock network of arbitrary size \cite{Daus20}.

Firstly, let us consider the iterative method, introduced in \cite{CH18}, to construct an Hermite interpolant starting from the basis of the renowned Floater-Hormann family of interpolants. This consists into considering primarily the basis functions $b_i(x)$ that are constructed, given $d\in\bbN$ and $x_0,\dots,x_n$ interpolation nodes, as 
$$ b_i(x) = \frac{w_i}{x-x_i}\bigg/ \sum_{j=0}^n \frac{w_j}{x-x_k}, i=0,\dots,n,$$
with 
$$ w_i = (-1)^{i+d} \sum_{j=\max(0,i-d)}^{\min(i,n-d)} \prod_{k=j, k\neq i}^{j+d} \frac{1}{|x_i-x_k|}, i=0,\dots,n.$$
Therefore, once are considered the functions
$$ b_{i,j}(x) = \frac{1}{j!} (x-x_i)^j b_i(x)^{j+1},$$
they can be used to construct an interpolant that satisfies the Hermite conditions.

The effectiveness of such construction for FH-interpolants has been investigated and in particular, in \cite{CH19}, it has been studied the Lebesgue constant at equidistant nodes and in \cite{CH20}, an estimate of the interpolation error has been given.

More specifically, we are going to firstly generalise this iterative Hermite construction and, later, use this method to introduce a barycentric rational trigonometric Hermite interpolant. 

The paper is construct as follows: in Section 2 is presented the general iterative approach. In Section 3 is analysed the Hermite interpolant constructed by using as basis function the one of Berrut's trigonometric interpolant and, moreover, is computed analytically the differentation matrix for the Hermite basis for general ordered nodes; this allows a fast construction and computation of the barycentric rational trigonometric Hermite interpolant. In Section 4 various numerical examples are presented and analysed. Finally, we conclude in Section 5.

\section{A general iterative approach}

Let us consider a set $\Omega \subset\bbR$ and some nodes $x_i\in \Omega$, $i=0,\dots,n$. Then, let us consider some basis functions $b_i$ construct by using the nodes $x_i$'s for an interpolant on $\Omega$ satisfying the Lagrange property $b_i(x_j) = \delta_{i,j}$.
We define, therefore, the following functions,  
$$
    b_{i,j}(x) :=\frac{1}{j!} d_i(x)^j b_i(x)^{j+1} 
$$
for a smooth function $d_i(x)$ that vanishes in $x_i$ and such that $d_i(x_j)\neq 0$, for $j\neq i$, and $d_i'(x_i)=1$. 

These functions will be the foundation to construct the Hermite interpolants, in fact we get the following.

\begin{lem} \label{LemGenBasis}
Let $\calX_n = \{ x_0,\dots,x_n\}\subset \Omega\subset\bbR$ be some ordered nodes and let $d_i(x)$ be a  function that vanishes in $x_j$ only if $j=i$ and such that $d_i'(x_i)=1$. Let $b_i(x)$, $i=0,\dots,n$ be a basis of some space $\calV$ that satisfy the Lagrange property at the nodes. Then, the functions $b_{i,j}$ defined as 
\begin{equation}\label{basis_hermite}
    b_{i,j}(x) :=\frac{1}{j!} d_i(x)^j b_i(x)^{j+1} 
\end{equation} 
satisfy
\begin{equation} \label{HermiteBasisCond}
    b_{i,j}^{(k)} (x_\ell) = \begin{system}
0, \quad\quad \text{ if }k<j, \\
\delta_{i,\ell}, \quad\, \text{ if }k=j.
\end{system}
\end{equation}
\end{lem}

\begin{proof}
For $j=0$ the statement follows directly from the Lagrange property of $b_i(x)$. Let us now consider $j>0$ and prove the statement by induction on $j$.

Let 
$$ c_i(x)= d_i(x) b_i(x); $$
then, it is clear that
$$ b_{i,j}(x) = \frac{1}{j} c_i(x) b_{i,j-1}(x). $$
Once we derive, we get by using the Leibniz rule that
$$ b_{i,j}^{(k)}(x) = \frac{1}{j} \sum_{s=0}^k \binom{k}{s} c_i(x)^{(k-s)}(x) b_{i,j-1}^{(s)}(x) .$$
Therefore, since $c_i(x_p) = 0 $ for $p=0,\dots,n$ and, since for the induction hypothesis we have that 
$$ b_{i,j-1}^{(s)}(x_p)=0 \text{ for } s<j-1,$$
we get 
$$ b_{i,j}^{(k)}(x_p) = \frac{1}{j} \sum_{s=j-1}^{k-1} \binom{k}{s} c_i^{(k-s)}(x_p) b_{i,j-1}^{(s)}(x_p). $$

Let us observe that the sum is empty in the case when $k<j$, whereas when $k=j$ we obtain
$$ b_{i,j}^{(j)}(x_p) = \frac{1}{j} \binom{j}{j-1} c_i'(x_p) b_{i,j-1}^{(j-1)}(x_p).$$

Finally, due to the induction hypothesis on $b_{i,j-1}^{(j-1)}(x_p)$, we notice that $b_{i,j}^{(j)}(x_p)=0$ if $p\neq i$ and when $p=i$ we have that the product is one since $c_i'(x_i)=1$ by construction. 
\end{proof}

\begin{remark}
Given a function $h_i(x)$ vanishing only in $x_i$ and with a non-zero derivative in $x_i$, we remark that we can always construct a function that satisfies the assumption of Lemma \ref{LemGenBasis} by setting 
$$ d_i(x)\coloneqq \frac{h_i(x)}{h_i'(x_i)}.$$
\end{remark}

Once we have the functions $b_{i,j}$ we can construct the Hermite interpolant starting from
$$ r_0(x) = \sum_{i=0}^n b_{i,0}(x) f_{i,0}, $$
which corresponds to the classical interpolant since $b_{i,0} = b_i$. For $j>0$, since we have that normally 
$$ r_{j-1}^{(j)}(x_\ell) \neq f_{\ell,j},$$
we introduce, at each step, a correction function $q_j$ which we add to the previous interpolants, in order to construct a new one that interpolates correctly the derivatives up to the $j$-th, i.e.,
\begin{equation} \label{HerInt}
    r_j(x) = r_{j-1}(x) + q_j(x), \quad j=1,\dots,m.
\end{equation} 
Therefore, we define the correction term as
\begin{equation}
    q_j(x) = \sum_{i=0}^n b_{i,j}(x) \Big( f_{i,j}-r^{(j)}_{j-1}(x_i) \Big), \quad j=1,\dots,m.
\end{equation}

The resulting function $r_m$, which we may write shortly as
$$ r_m(x) \coloneqq \sum_{j=0}^m \sum_{i=0}^n b_{i,j}(x) g_{i,j}$$
for proper $g_{i,j}$, construct via the iterative approach described above, thanks to Lemma \ref{LemGenBasis}, will satisfy the Hermite conditions. In fact, we have the following theorem.

\begin{thm}
    Let $\calX_n = \{ x_0,\dots,x_n\}\subset \Omega\subset\bbR$ be some ordered nodes and let $d_i(x)$ be a  function that vanishes in $x_j$ only if $j=i$ and such that $d_i'(x_i)=1$. Let $b_i(x)$, $i=0,\dots,n$ be a basis of some space $\calV$ that satisfy the Lagrange property at the nodes and let $b_{i,j}$ be defined as \eqref{basis_hermite}. Then, the function $r_m$, defined as
    \begin{equation}
        r_m(x) \coloneqq \sum_{j=0}^m \sum_{i=0}^n b_{i,j}(x) g_{i,j},
    \end{equation}
    with 
    \begin{equation}
        g_{i,j} = \begin{system}
            f_{i,0} \qquad\qquad\qquad \text{ if }j=0 \\
            f_{i,j}-r^{(j)}_{j-1}(x_i) \quad\; \text{ if }j>0,
        \end{system}
    \end{equation}
satisfies the Hermite conditions \eqref{HermCond}
\end{thm}

\begin{proof}
It is straightforward that, $$ r_m^{(k)}(x_\ell) = \sum_{i=0}^n \sum_{j=0}^m b_{ i,j}^{(k)}(x)_{\lvert_{x=x_\ell}} g_{i,j}, $$
which due to Lemma \ref{LemGenBasis} is 
\begin{align*}
    r_m^{(k)}(x_\ell) &= \sum_{i=0}^n \sum_{j=0}^k b_{ i,j}^{(k)}(x)_{\lvert_{x=x_\ell}} g_{i,j} \\
    &= g_{\ell,k} + \sum_{i=0}^n \sum_{j=0}^{k-1} b_{ i,j}^{(k)}(x)_{\lvert_{x=x_\ell}} g_{i,j} \\
    &= g_{\ell,k} + r_{k-1}^{(k)}(x_\ell) \\
    &= f_{\ell,k} .
\end{align*}
\end{proof}

\section{A barycentric rational trigonometric interpolant} \label{sec3}

Barycentric rational trigonometric interpolants are suitable for approximating periodic functions and many of them have been introduced in the last years. 

The barycentric form of the classical interpolant at equidistant nodes has been introduced by Henrici \cite{Henrici79}, and, later, Berrut \cite{Berrut88} proposed to use the same function as the Henrici's interpolant but for other ordered nodes.

More recently, Bandiziol and De Marchi \cite{BD19} proposed a barycentric rational trigonometric interpolant constructed similarly as the well-known Floater-Hormann interpolant. Moreover, Baddoo \cite{Baddoo21} introduced a trigonometric equivalent of the renowned algorithm AAA (see \cite{Nakatsukasa18}), the AAAtrig, which construct a barycentric trigonometric rational approximant selecting the nodes progressively via a greedy algorithm.

In this section, we focus into Berrut's trigonometric interpolant \cite{Berrut88}. Let us consider $n$ ordered nodes $0\leq \theta_0<\dots<\theta_{n-1}<2\pi$, then, we define the interpolant as 
\begin{equation} \label{BaryTrigEqui}
    T_n(\theta) = \frac{ \sum_{i=0}^{n-1} (-1)^i \cst\left( \frac{\theta-\theta_i}{2}  \right) f(\theta_i)}{\sum_{i=0}^{n-1} (-1)^i \cst\left( \frac{\theta-\theta_i}{2}  \right)},
\end{equation}
where the function $\cst$ is 
\begin{equation} \label{cst}
    \cst(\theta)= \begin{system} \csc(\theta), \quad \text{ if } n \text{ is odd,} \\
\ctg(\theta), \quad \text{ if } n \text{ is even}.
\end{system}
\end{equation}
It has no poles in $[0,2\pi)$ and enjoys a logarithmic growth of its Lebesgue constant for a wide family of interpolation nodes \cite{BE21}.
This class of nodes includes also the images of equidistant nodes via a conformal map, such as those in \cite{Baddoo19,BE20}.
As said earlier, it corresponds to the classical trigonometric interpolant when the nodes are equidistant \cite{Henrici79} and, therefore, it enjoys its properties, such as the exponential convergence. In addition, the interpolant converges exponentially also when the nodes are images of equidistant nodes via a conformal map \cite{Baltensperger02} and it has been used effectively in \cite{Berrut22} to interpolate functions on two-dimensional starlike domains.

In order to construct the Hermite interpolant and to retain the periodic behaviour, we may choose 
\begin{equation}
    d_i(\theta) = 2 \sin\left( \frac{\theta-\theta_i}{2}\right), 
\end{equation}
which clearly satisfies the conditions we need to define the functions $b_{i,j}$. 

In particular, we get
\begin{align}
    b_{i,j}(\theta) &\coloneqq \frac{2^j}{j!} \sin^j\left( \frac{\theta-\theta_i}{2}\right)(-1)^{ij+i}\frac{ \cst^{j+1} \left( \frac{\theta-\theta_i}{2}\right)}{\left( \sum_{k=0}^{n-1} (-1)^k \cst \left( \frac{\theta-\theta_k}{2}\right)\right)^{j+1}} \nonumber \\
    &= \frac{2^j}{j!} (-1)^{ij+i}\frac{ \cst \left( \frac{\theta-\theta_i}{2}\right)}{\left( \sum_{k=0}^{n-1} (-1)^k \cst \left( \frac{\theta-\theta_k}{2}\right)\right)^{j+1}} \cdot \begin{system}
1,   \qquad\qquad\;\; \text{ if } n \text{ is odd,}  \\
\cos^{j}\left( \frac{\theta-\theta_i}{2}\right),    \text{ if } n \text{ is even,}
    \end{system} \nonumber \\
    & = \frac{2^j (-1)^{i(j+1)}}{j!}  \cst\bigg(\frac{\theta -\theta_i}{2}\bigg)   \begin{system}
\frac{1}{ \left( \sum_{k} \frac{(-1)^k}{\sin\left( (\theta-\theta_k)/2 \right)} \right)^{j+1}}   \qquad \text{ if } n \text{ is odd,} \\
\\
 \frac{  \cos^j\big(\frac{\theta-\theta_i}{2}\big)  }{ \left( \sum_{k} \frac{(-1)^k}{\sin\left( (\theta-\theta_k)/2 \right)} \right)^{j+1} }  \quad\quad   \text{ if } n \text{ is even}, 
    \end{system} \nonumber \\
& = \dfrac{\frac{2^j (-1)^{i(j+1)}}{j!}  \cst\bigg(\frac{\theta -\theta_i}{2}\bigg) \, \eta_{i,j}^{(n)}(\theta),}{\left( \sum_{k=0}^{n-1} \frac{(-1)^k}{\sin\left( (\theta-\theta_k)/2 \right)} \right)^{j+1}} \label{HermiteTrigoBasis}
\end{align}
with 
$$\eta_{i,j}^{(n)}(\theta) = \begin{system}
    1\qquad \qquad\quad\,\, \text{ if $n$ is odd,} \\
     \cos^j\big(\frac{\theta-\theta_i}{2}\big) \quad \text{ if $n$ is even.}
\end{system}$$

This formulation will allow us to calculate analytically the differential matrix of the Hermite interpolant 
\begin{equation}\label{TrigHermite}
    t_m(\theta) = \sum_{i=0}^{n-1} \sum_{j=0}^m b_{i,j}(\theta) g_{i,j},
\end{equation} 
which will be useful for a fast computation of the interpolant.

\begin{thm}
Let us consider the basis functions \eqref{HermiteTrigoBasis}, then  when $s>j$ the differential matrix, whose elements are $(D^{s}_j)_{ik} = b_{k,j}^{(s)}(\theta_i)$, is defined as
\begin{align} 
      {(D_{j}^{s})}_{ik} = & \frac{1}{\sin\big(\frac{\theta_i-\theta_k}{2}\big)}\bigg((-1)^{(k-i)(j+1)} \sum_{q=0}^{\lfloor \frac{s-j-1}{2} \rfloor} \binom{s}{2q+1} \frac{(-1)^q}{2^{2q+1}} (D_{j}^{s-2q-1})_{ii} \nonumber \\
      & - \sum_{q=j}^{s-1} \binom{s}{s-q} \bigg(\sin\bigg(\frac{\theta-\theta_k}{2}\bigg)\bigg)^{(s-q)}_{|_{\theta=\theta_i}} (D_{j}^{q})_{ik} \bigg) 
\end{align}
when $n$ is odd and
\begin{align}
    {(D_{j}^{s})}_{ik}  = & \frac{1}{\tan\big(\frac{\theta_i-\theta_k}{2}\big)}\bigg((-1)^{(k-i)(j+1)} \sum_{q=0}^{\lfloor \frac{s-j}{2}\rfloor} \binom{s}{2q+1} 2^{-2q-1} \cdot \nonumber \\
    & \quad \cdot \bigg(\sum_{j=0}^{q} a(2q+1,2j+1) \bigg) (D_{j}^{s-2q-1})_{ii}  \nonumber \\ 
    & \quad - \sum_{q=j}^{s-1} \binom{s}{s-q} \bigg(\tan\bigg(\frac{\theta-\theta_k}{2}\bigg)\bigg)^{(s-q)}_{|_{\theta=\theta_i}} (D_{j}^{q})_{ik} \bigg).
\end{align}
when $n$ is even, with
$$ a(p,q) = (-1)^{\frac{1}{2}(p-\frac{3+(-1)^p}{2})}\,2\sum_{j=0}^{\frac{p-q-1}{2}} (-1)^{\frac{p-q-1}{2}-j} \binom{p+1}{j} \bigg( \frac{p-q-1}{2}-j+1\bigg)^p. $$
\end{thm}

\begin{proof}
Therefore, we can write the basis as
\begin{equation} \label{TrigHermiteBas}
    b_{k,j}(\theta) := u_{k,j} \alpha(\theta) \cst\Big(\frac{\theta-\theta_k}{2}\Big)
\end{equation}
where 
\begin{equation}\label{baryweights}
    u_{k,j} \coloneqq \frac{2^j (-1)^{k(j+1)}}{j!} 
\end{equation}
and, if $n$ is odd,
$$ \alpha(\theta) := \frac{1}{\Big(\sum_{i=0}^{n-1} (-1)^i \cst((\theta-\theta_i)/2)\Big)^{j+1} } $$
and, if $n$ is even,
$$ \alpha(\theta) :=  \frac{\cos^j(\frac{\theta-\theta_k}{2})}{ \Big( \sum_{i=0}^{n-1} (-1)^i \cst\left( (\theta-\theta_i)/2 \right) \Big)^{j+1} }  .$$

In this way we can compute the differentation matrices as done for the classical barycentric trigonometric interpolant in \cite{Baltensperger02}, which will be useful to compute the values of the derivatives of the interpolant $r_m$ at the nodes in the previous iterations.

If $n$ is odd, considering equation \eqref{TrigHermiteBas}, differentiating both side and evaluating at the points $x_i$, we get that
\begin{align}
    \big(\alpha(\theta)\big)^{(s)}_{|_{\theta=\theta_k}} &= \frac{1}{u_{k,j}} \sum_{q=0}^s \binom{s}{q} \bigg(\sin\bigg(\frac{\theta-\theta_k}{2}\bigg)\bigg)^{(q)}\big|_{\theta=\theta_k} \big(b_{k,j}(\theta)\big)^{(s-q)}\big|_{\theta=\theta_k} \nonumber \\
    &= \frac{1}{u_{k,j}} \sum_{q=0}^s \binom{s}{q} \bigg(\sin\bigg(\frac{\theta-\theta_k}{2}\bigg)\bigg)^{(q)}\big|_{\theta=\theta_k} (D_j^{s-q})_{kk} \nonumber \\
    &= \frac{1}{u_{k,j}} \sum_{q=1}^s \binom{s}{q} \bigg(\sin\bigg(\frac{\theta-\theta_k}{2}\bigg)\bigg)^{(q)}\big|_{\theta=\theta_k} (D_{j}^{s-q})_{kk}, \label{DifMat1_odd}
\end{align}
where we defined $(D^{s}_j)_{ik} = b_{k,j}^{(s)}(\theta_i)$.\\
Moreover, by evaluating the derivative in different nodes we have
\begin{align}
    \big(\alpha(\theta)\big)^{(s)}\big|_{\theta=\theta_i} &= \frac{1}{u_{k,j}} \sum_{q=0}^s \binom{s}{q} \bigg(\sin\bigg(\frac{\theta-\theta_k}{2}\bigg)\bigg)^{(q)}\big|_{\theta=\theta_i} \big(b_{k,j}(\theta)\big)^{(s-q)}\big|_{\theta=\theta_i} \nonumber \\ 
    &= \frac{1}{u_{k,j}} \sum_{q=0}^s \binom{s}{q} \bigg(\sin\bigg(\frac{\theta-\theta_k}{2}\bigg)\bigg)^{(q)}\big|_{\theta=\theta_i} (D_{j}^{s-q})_{ik}. \label{DifMat2_odd}
\end{align}
In \eqref{DifMat2_odd} we can isolate the term $q=0$ of the sum, that gives,
\begin{align}
    \big(\alpha(\theta)\big)^{(s)}\big|_{\theta=\theta_i} =& \frac{1}{u_{k,j}} \sum_{q=1}^s \binom{s}{q} \bigg(\sin\bigg(\frac{\theta-\theta_k}{2}\bigg)\bigg)^{(q)}\big|_{\theta=\theta_i} (D_{j}^{s-q})_{ik} \nonumber \\
    &+ \frac{1}{u_{k,j}} \bigg(\sin\bigg(\frac{\theta_i-\theta_k}{2}\bigg)\bigg) (D_{j}^{s})_{ik},\label{DifMat3}
\end{align}
which together with \eqref{DifMat1_odd} gives us the following recursive formula for the differentation matrix
\begin{align}
    {(D_{j}^{s})}_{ik} & = \frac{1}{\sin\big(\frac{\theta_i-\theta_k}{2}\big)}\bigg( u_{k,j} \big(\alpha(\theta)\big)^{(s)}\big|_{\theta=\theta_i} \nonumber \\
    & \qquad\qquad\qquad\quad - \sum_{q=1}^s \binom{s}{q} \bigg(\sin\bigg(\frac{\theta-\theta_k}{2}\bigg)\bigg)^{(q)}\big|_{\theta=\theta_i} (D_{j}^{s-q})_{ik} \bigg) \\ 
    &= \frac{1}{\sin\big(\frac{\theta_i-\theta_k}{2}\big)}\bigg( \frac{u_{k,j}}{u_{i,j}} \sum_{q=1}^s \binom{s}{q} \bigg(\sin\bigg(\frac{\theta-\theta_i}{2}\bigg)\bigg)^{(q)}\big|_{\theta=\theta_i} (D_{j}^{s-q})_{ii} \nonumber \\
    &\qquad\qquad\quad \qquad- \sum_{q=1}^s \binom{s}{q} \bigg(\sin\bigg(\frac{\theta-\theta_k}{2}\bigg)\bigg)^{(q)}\big|_{\theta=\theta_i} (D_{j}^{s-q})_{ik} \bigg). \label{DifMat4_odd}
\end{align}
Notice that since the terms $b_{i,j}$ satisfy \eqref{HermiteBasisCond}, $D_j^j$ corresponds to the identity matrix and $D_{j}^{s}=0$ when $s<j$. Then, the formula becomes
\begin{align} 
      {(D_{j}^{s})}_{ik} = & \frac{1}{\sin\big(\frac{\theta_i-\theta_k}{2}\big)}\bigg(\frac{u_{k,j}}{u_{i,j}} \sum_{q=j}^{s-1} \binom{s}{s-q} \bigg(\sin\bigg(\frac{\theta-\theta_i}{2}\bigg)\bigg)^{(s-q)}\big|_{\theta=\theta_i} (D_{j}^{q})_{ii} \nonumber \\
      & - \sum_{q=j}^{s-1} \binom{s}{s-q} \bigg(\sin\bigg(\frac{\theta-\theta_k}{2}\bigg)\bigg)^{(s-q)}\big|_{\theta=\theta_i} (D_{j}^{q})_{ik} \bigg). \label{RecFormu_even}
\end{align}

Similarly, when $n$ is even, we have that 
\begin{align}
    \big(\alpha(\theta)\big)^{(s)}\big|_{\theta=\theta_k} &= \frac{1}{u_{k,j}} \sum_{q=1}^s \binom{s}{q} \bigg(\tan\bigg(\frac{\theta-\theta_k}{2}\bigg)\bigg)^{(q)}\big|_{\theta=\theta_k} (D_{j}^{s-q})_{kk} \label{DifMat1_even} \\
    \big(\alpha(\theta)\big)^{(s)}\big|_{\theta=\theta_i} &= \frac{1}{u_{k,j}} \sum_{q=0}^s \binom{s}{q} \bigg(\tan\bigg(\frac{\theta-\theta_k}{2}\bigg)\bigg)^{(q)}\big|_{\theta=\theta_i} (D_{j}^{s-q})_{ik}. \label{DifMat2_even}
\end{align}
Hence, 
\begin{align}
    {(D_{j}^{s})}_{ik} = & \frac{1}{\tan\big(\frac{\theta_i-\theta_k}{2}\big)}\bigg(\frac{u_{k,j}}{u_{i,j}} \sum_{q=j}^{s-1} \binom{s}{s-q} \bigg(\tan\bigg(\frac{\theta-\theta_i}{2}\bigg)\bigg)^{(s-q)}\big|_{\theta=\theta_i} (D_{j}^{q})_{ii} \nonumber \\
      & - \sum_{q=j}^{s-1} \binom{s}{s-q} \bigg(\tan\bigg(\frac{\theta-\theta_k}{2}\bigg)\bigg)^{(s-q)}\big|_{\theta=\theta_i} (D_{j}^{q})_{ik} \bigg). \label{RecFormu_odd}
\end{align}

Moreover, since
$$ \bigg(\sin\bigg(\frac{\theta-\theta_k}{2}\bigg)\bigg)^{(q)} = \begin{system}
\frac{(-1)^{\frac{q}{2}}}{2^q} \sin\big(\frac{\theta-\theta_k}{2}\big), \qquad\quad \text{ if }q \text{ even} \\
\frac{(-1)^{\frac{q-1}{2}}}{2^q} \cos\big(\frac{\theta-\theta_k}{2}\big), \qquad \text{ if }q \text{ odd} 
\end{system}$$
and since the $j$-th derivative of $\tan((\theta-\theta_i)/2)$ computed in $\theta_i$ is a given by (see \cite{Qi15})
$$ \bigg( (\tan\bigg(\frac{\theta-\theta_i}{2}\bigg)\bigg)^{(n)}\big|_{\theta=\theta_i} = \begin{system}
0, \qquad\qquad\qquad\qquad\qquad\qquad\;  \text{ if }n = 2k, \\
2^{-n} \sum_{j=0}^{k-1} a(2k-1,2j+1), \quad \text{ if } n = 2k - 1 ,
\end{system}$$
for some $k\in\bbN$ and with 
$$ a(p,q) = (-1)^{\frac{1}{2}(p-\frac{3+(-1)^p}{2})}\,2\sum_{j=0}^{\frac{p-q-1}{2}} (-1)^{\frac{p-q-1}{2}-j} \binom{p+1}{j} \bigg( \frac{p-q-1}{2}-j+1\bigg)^p $$
we can simplify equation \eqref{RecFormu_even} to
\begin{align} 
      {(D_{j}^{s})}_{ik} = & \frac{1}{\sin\big(\frac{\theta_i-\theta_k}{2}\big)}\bigg((-1)^{(k-i)(j+1)} \sum_{q=0}^{\lfloor \frac{s-j-1}{2} \rfloor} \binom{s}{2q+1} \frac{(-1)^q}{2^{2q+1}} (D_{j}^{s-2q-1})_{ii} \nonumber \\
      & - \sum_{q=j}^{s-1} \binom{s}{s-q} \bigg(\sin\bigg(\frac{\theta-\theta_k}{2}\bigg)\bigg)^{(s-q)}_{|_{\theta=\theta_i}} (D_{j}^{q})_{ik} \bigg) 
\end{align}
and equation \eqref{RecFormu_odd} to
\begin{align}
    {(D_{j}^{s})}_{ik}  = & \frac{1}{\tan\big(\frac{\theta_i-\theta_k}{2}\big)}\bigg((-1)^{(k-i)(j+1)} \sum_{q=0}^{\lfloor \frac{s-j}{2}\rfloor} \binom{s}{2q+1} 2^{-2q-1} \cdot \nonumber \\
    & \quad \cdot \bigg(\sum_{j=0}^{q} a(2q+1,2j+1) \bigg) (D_{j}^{s-2q-1})_{ii}  \nonumber \\ 
    & \quad - \sum_{q=j}^{s-1} \binom{s}{s-q} \bigg(\tan\bigg(\frac{\theta-\theta_k}{2}\bigg)\bigg)^{(s-q)}_{|_{\theta=\theta_i}} (D_{j}^{q})_{ik} \bigg).
\end{align}
\end{proof}

The interpolant \eqref{TrigHermite} then can be written in the following extended form
\begin{align}
    t_m(\theta) &= \sum_{i=0}^{n-1} \sum_{j=0}^m b_{i,j}(\theta) \bigg( f_{i,j} - \sum_{s=0}^{j-1} \sum_{k=0}^{n-1} (D_s^j)_{ik} f_{k,s} \bigg)
\end{align} 
or the following form linear on the data 
\begin{align}
    t_m(\theta) &= \sum_{i=0}^{n-1}  \sum_{j=0}^m \bigg( b_{i,j}(\theta) - \sum_{s=j+1}^{m} \sum_{k=0}^{n-1} b_{k,s}(\theta) (D_j^s)_{ki} \bigg) f_{i,j},
\end{align} 
which is specifically 
\begin{align} 
    t_m(\theta) =& \Bigg(\sum_{i=0}^{n-1}  \sum_{j=0}^m \bigg( u_{i,j} \cst\Big(\frac{\theta-\theta_i}{2}\Big) \eta_{i,j}^{(n)}(\theta) \left(\sum_{i=0}^{n-1} \frac{(-1)^i}{\sin\left( \frac{\theta-\theta_i}{2} \right)}\right)^{m-j} \nonumber \\
    & - \sum_{s=j+1}^{m} \sum_{k=0}^{n-1} u_{k,s}\cst\Big(\frac{\theta-\theta_k}{2}\Big) \eta_{i,j}^{(n)}(\theta) \left(\sum_{i=0}^{n-1} \frac{(-1)^i}{\sin\left( \frac{\theta-\theta_i}{2} \right)}\right)^{m-s} (D_j^s)_{ki} \bigg) f_{i,j} \Bigg) \nonumber \\
    & \Bigg/  \left( \sum_{i=0}^{n-1} \frac{(-1)^i}{\sin\left( \frac{\theta-\theta_i}{2} \right)} \right)^{m+1},
\end{align} 
with $u_{i,j}$ defined as in \eqref{baryweights}.

Once we consider the interpolant of the function $f\equiv 1$, the interpolant resolves to
\begin{align} 
    t_m(\theta) =& \Bigg(\sum_{i=0}^{n-1}  \sum_{j=0}^m \bigg( \sum_{k=0}^{n-1} u_{k,j} (D_0^j)_{ki} \bigg) \cst\Big(\frac{\theta-\theta_i}{2}\Big) \left(\sum_{i=0}^{n-1} \frac{(-1)^i}{\sin\left( \frac{\theta-\theta_i}{2} \right)}\right)^{m-j} \Bigg) \nonumber \\
    & \Bigg/  \left( \sum_{i=0}^{n-1} \frac{(-1)^i}{\sin\left( \frac{\theta-\theta_i}{2} \right)} \right)^{m+1},
\end{align} 
since the derivatives of $f$ are zeros, and $D_0^0$ is the identity matrix.

Therefore we get the barycentric form of the interpolant, that is
\begin{align} 
    t_m(\theta) =& \Bigg(\sum_{i=0}^{n-1}  \sum_{j=0}^m \bigg( u_{i,j} \cst\Big(\frac{\theta-\theta_i}{2}\Big) \eta_{i,j}^{(n)}(\theta) \left(\sum_{i=0}^{n-1} \frac{(-1)^i}{\sin\left( \frac{\theta-\theta_i}{2} \right)}\right)^{m-j} \nonumber \\
    & - \sum_{s=j+1}^{m} \sum_{k=0}^{n-1} u_{k,s}\cst\Big(\frac{\theta-\theta_k}{2}\Big) \eta_{i,j}^{(n)}(\theta) \left(\sum_{i=0}^{n-1} \frac{(-1)^i}{\sin\left( \frac{\theta-\theta_i}{2} \right)}\right)^{m-s} (D_j^s)_{ki} \bigg) f_{i,j} \Bigg) \nonumber \\
    & \Bigg/  \Bigg(\sum_{i=0}^{n-1}  \sum_{j=0}^m \bigg( \sum_{k=0}^{n-1} u_{k,j} (D_0^j)_{ki} \bigg) \cst\Big(\frac{\theta-\theta_i}{2}\Big) \left(\sum_{i=0}^{n-1} \frac{(-1)^i}{\sin\left( \frac{\theta-\theta_i}{2} \right)}\right)^{m-j} \Bigg)
\end{align}

In addition, we remark that, when we compute the matrix for $s=j+1$, the elements simplify to the following
\begin{align} 
      {(D_{j}^{j+1})}_{ik} = & \frac{(-1)^{(j+1)(k-i)} \, (j+1)}{2}\cst\bigg(\frac{\theta_i-\theta_k}{2}\bigg), \qquad k\neq i \label{RecFormu2_Jp1_odd}, 
\end{align}
which for $j=0$ corresponds to the values in \cite{Baltensperger02}; for the diagonal elements, we use the relation $\sum_{i=0}^{n-1} b_{i,j}^{(j)}(\theta) = 1$, since we want to reproduce the constants for the $j$-th derivate. This implies that 
$$ \sum_{i=0}^{n-1} b_{i,j}^{(j+1)}(\theta) = 0, $$
so that we obtain 
$$ (D_{j}^{j+1})_{ii} = - \sum_{\substack{k=0 \\ k\neq i}}^{n-1} (D_{j}^{j+1})_{ik}. $$ 

Moreover, since for the construction of the Hermite interpolant \eqref{TrigHermite}, we need to compute $r^{(s)}_{j-1}(x_i)$ and therefore $D_{j}^{s}$, in the next section, we decided to proceed by computing the differentiation matrices of higher order as 
\begin{equation}
    D_j^{s} = (D_j^{j+1})^{s-j} \qquad s\geq j+1.
\end{equation} 
as done for the polynomial case in \cite{Welfert97}.

\section{Numerical Experiments}

Here we present some numerical tests, some Matlab demos to reproduce the experiments are available at \\
\texttt{https://github.com/gelefant/TrigonometricHermite} 

We considered the following periodic test functions in $[0,2\pi)$,
\begin{align*}
    f_1(\theta) &= e^{2\sin(\theta)+\cos(\theta)}, \\
    f_2(\theta) &= \cos(3\theta) + \log(\cos(\theta)+1.5),
\end{align*} 
and we interpolated them by means of the interpolant \eqref{TrigHermite} with the basis presented in Section \ref{sec3} with $N$ equidistant nodes.

\begin{figure}
    \centering
    \graphicspath{ {Images/}}
    \subfigure[]
        {\includegraphics[height=1.84in]{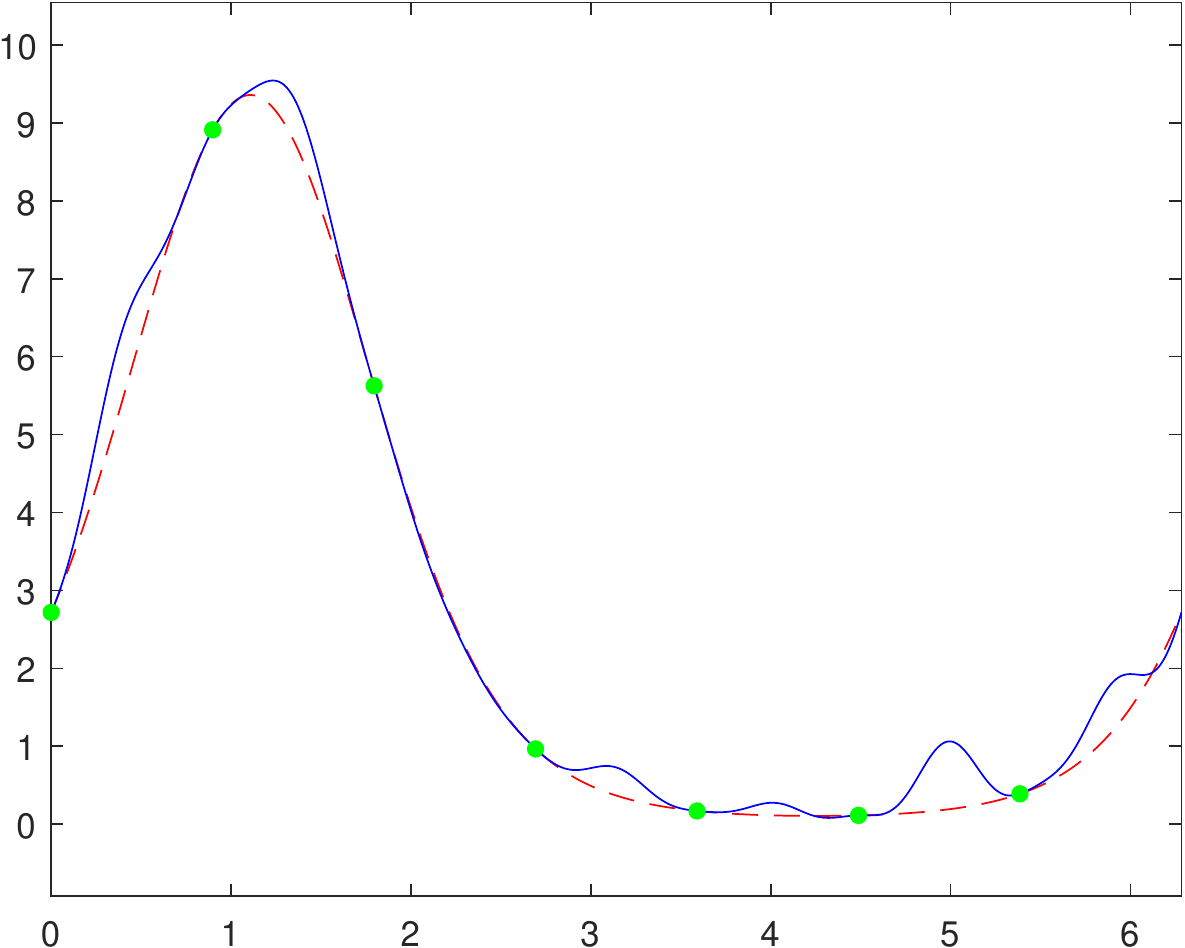}}
     ~ 
    \subfigure[]
        {\includegraphics[height=1.84in]{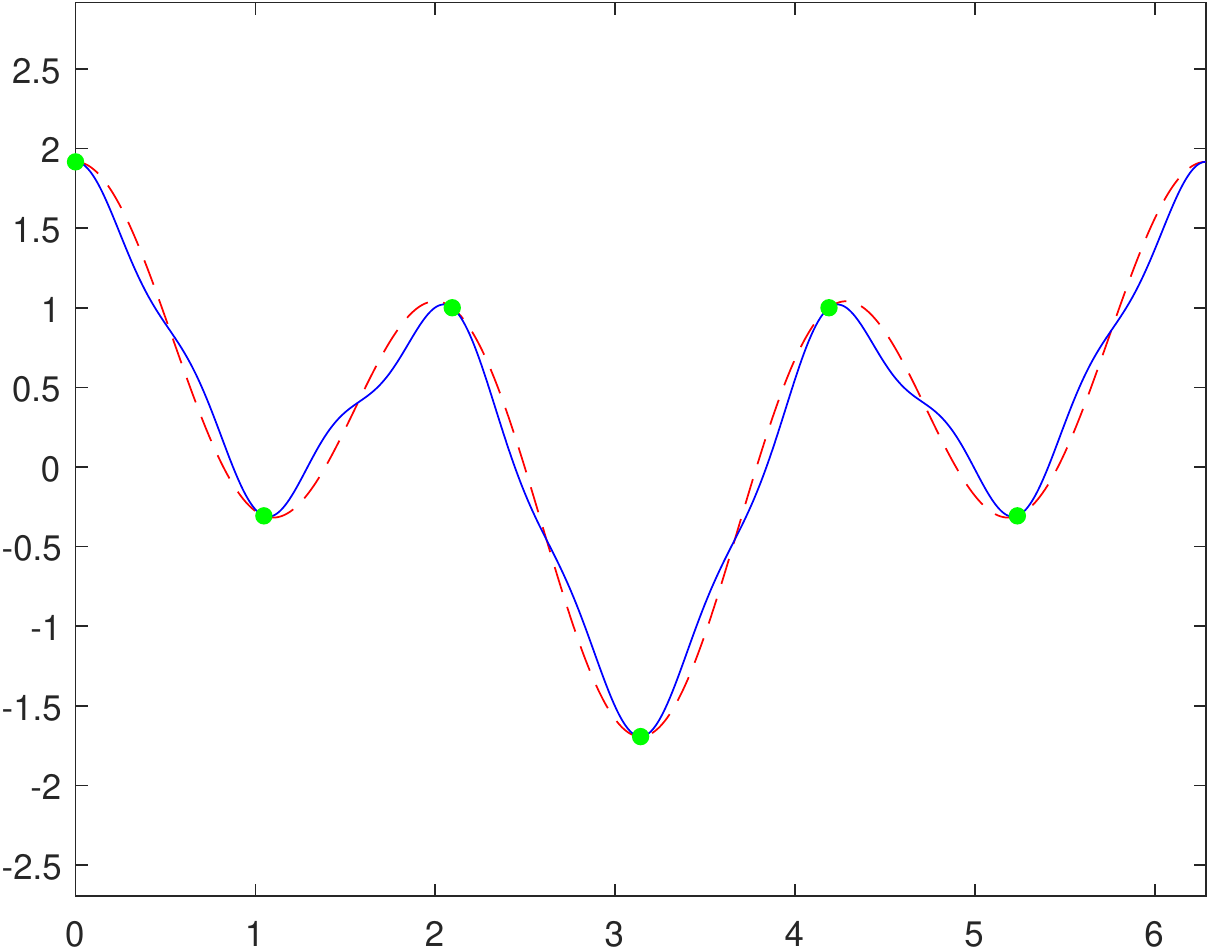}}
    \caption{The functions $f_1$ (left) and $f_2$ (right) in red dashed and in blue line the interpolant $t_3$ at $7$ equidistant nodes (left) and the interpolant $t_2$ at $6$ equidistant nodes (right). Interpolation nodes are labeled in green.   }
    \label{Fig1}
\end{figure}

We tested the interpolants which satisfy the Hermite conditions from one to four derivatives. In Figure \ref{Fig1} we display both test functions together with two interpolants, one satisfying the Hermite conditions up to the second derivative and one satisfying the conditions up to the third derivative.

\begin{figure}
    \centering
    \graphicspath{ {Images/}}
    \subfigure[]
        {\includegraphics[height=1.80in]{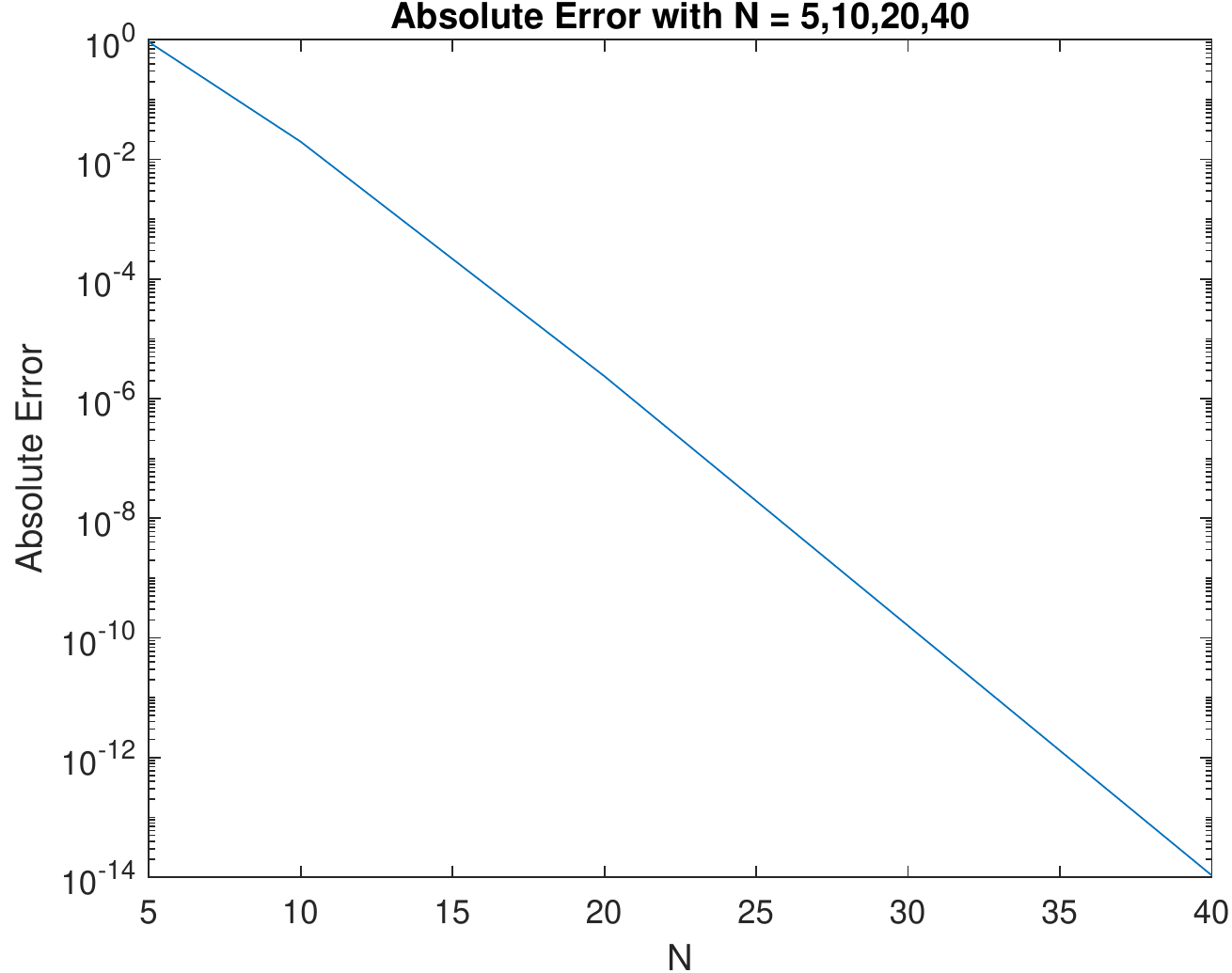}}
     ~ 
    \subfigure[]
        {\includegraphics[height=1.80in]{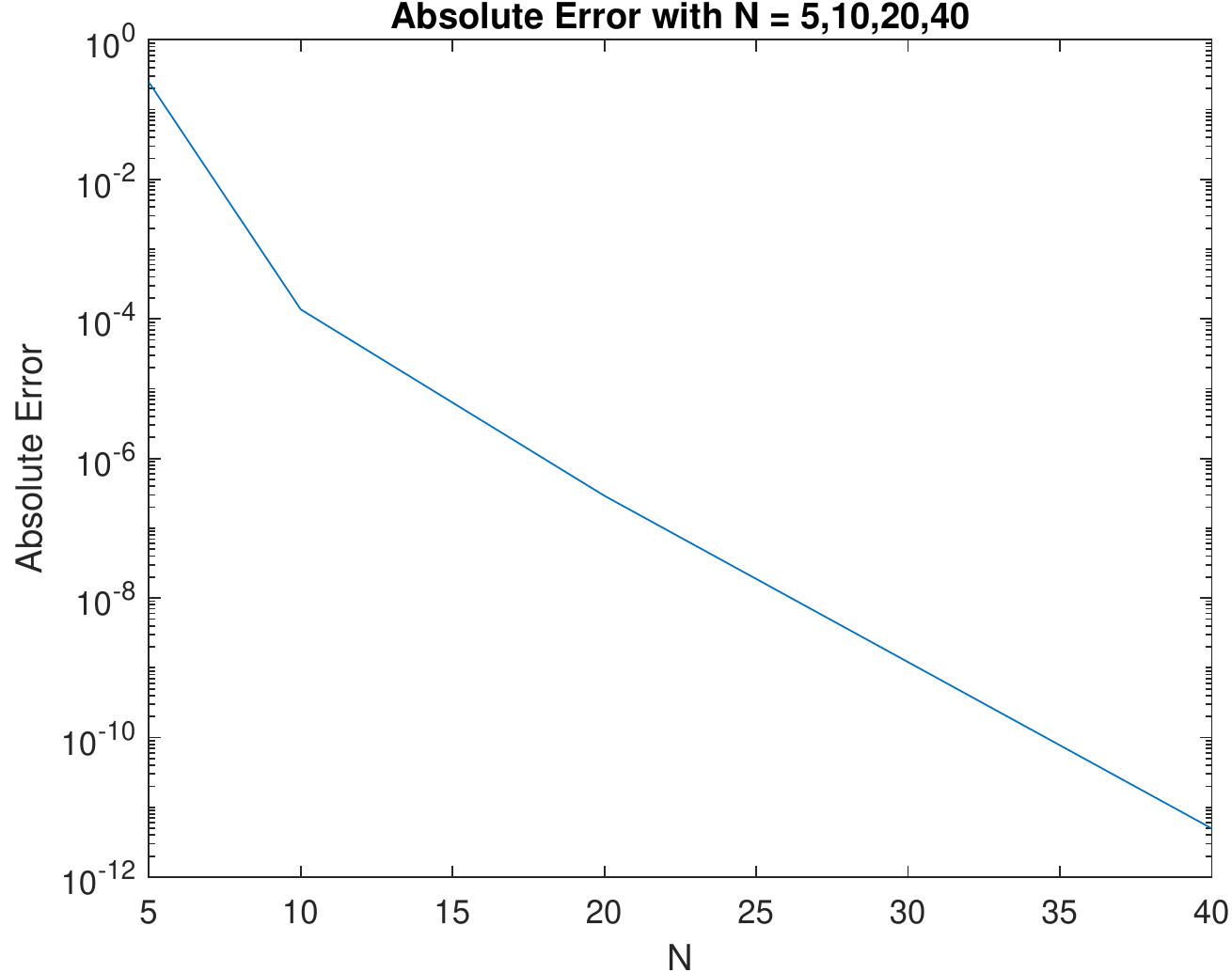}}
    \caption{Convergence of $t_1$ for $f_1$ (left) and $f_2$ (right) and $N=5,10,20,40$.  }
    \label{Fig_Conv_H1}
\end{figure}

\begin{figure}
    \centering
    \graphicspath{ {Images/}}
    \subfigure[]
        {\includegraphics[height=1.80in]{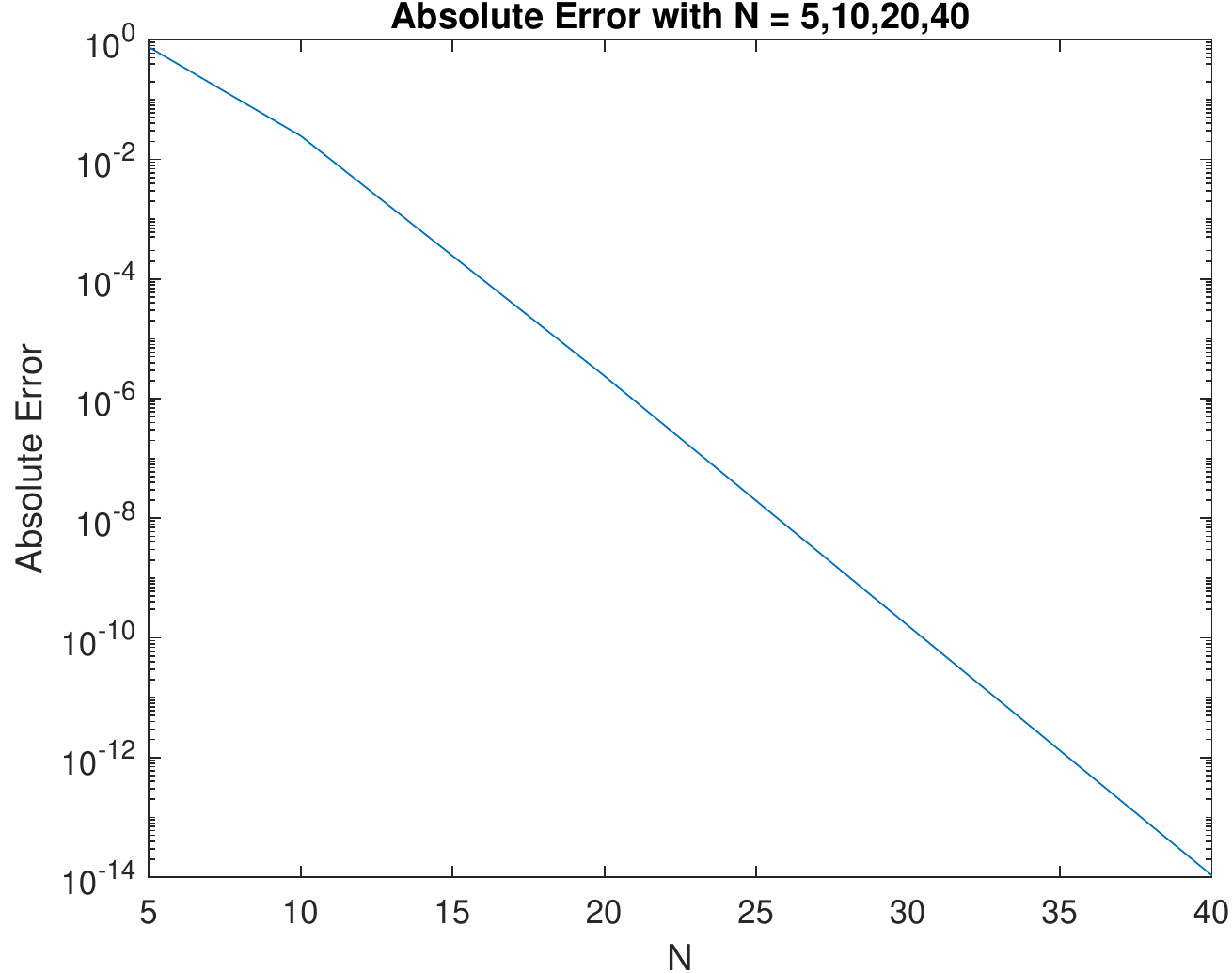}}
     ~ 
    \subfigure[]
        {\includegraphics[height=1.80in]{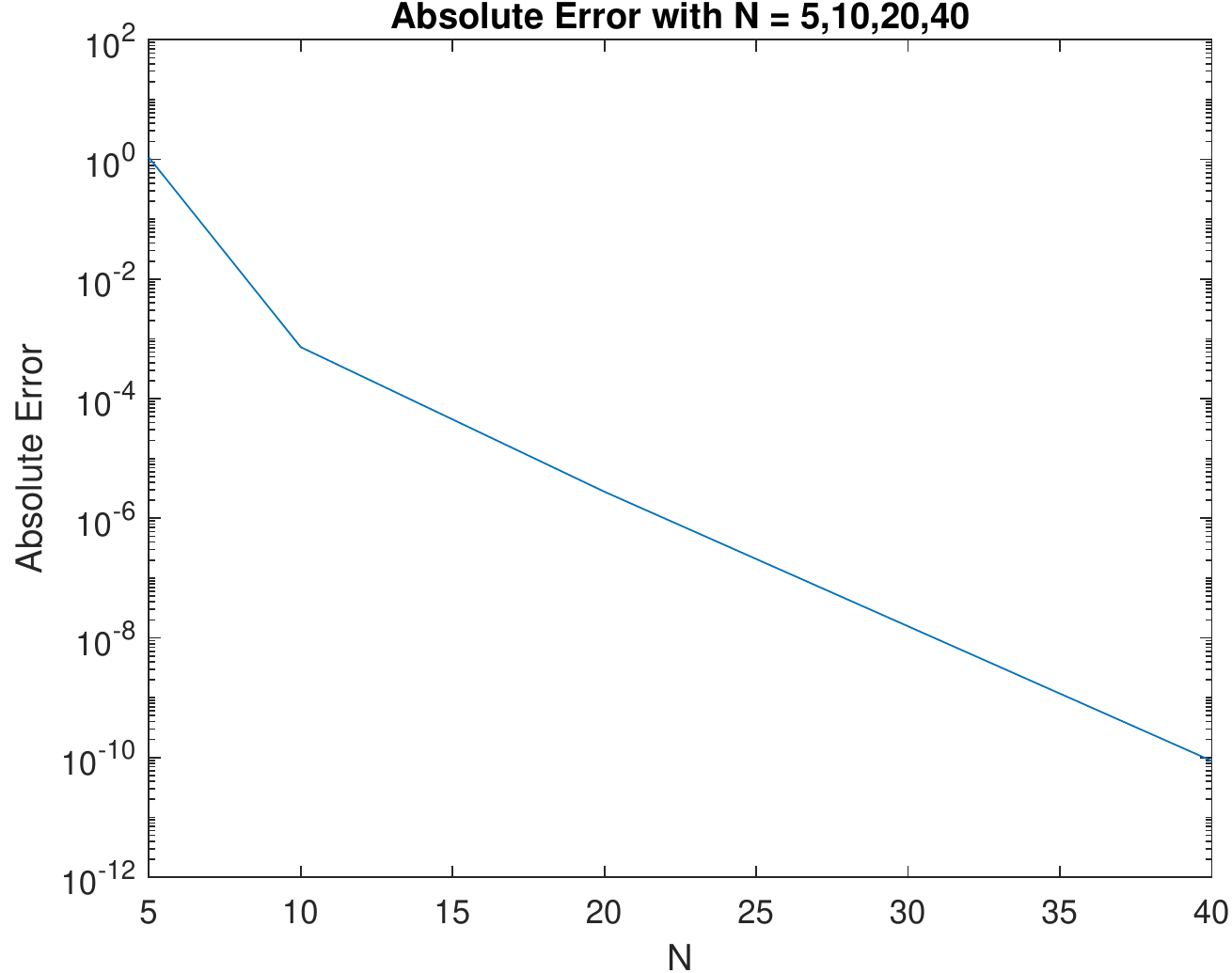}}
    \caption{Convergence of $t_2$ for $f_1$ (left) and $f_2$ (right) and $N=5,10,20,40$.  }
    \label{Fig_Conv_H2}
\end{figure}

\begin{figure}
    \centering
    \graphicspath{ {Images/}}
    \subfigure[]
        {\includegraphics[height=1.80in]{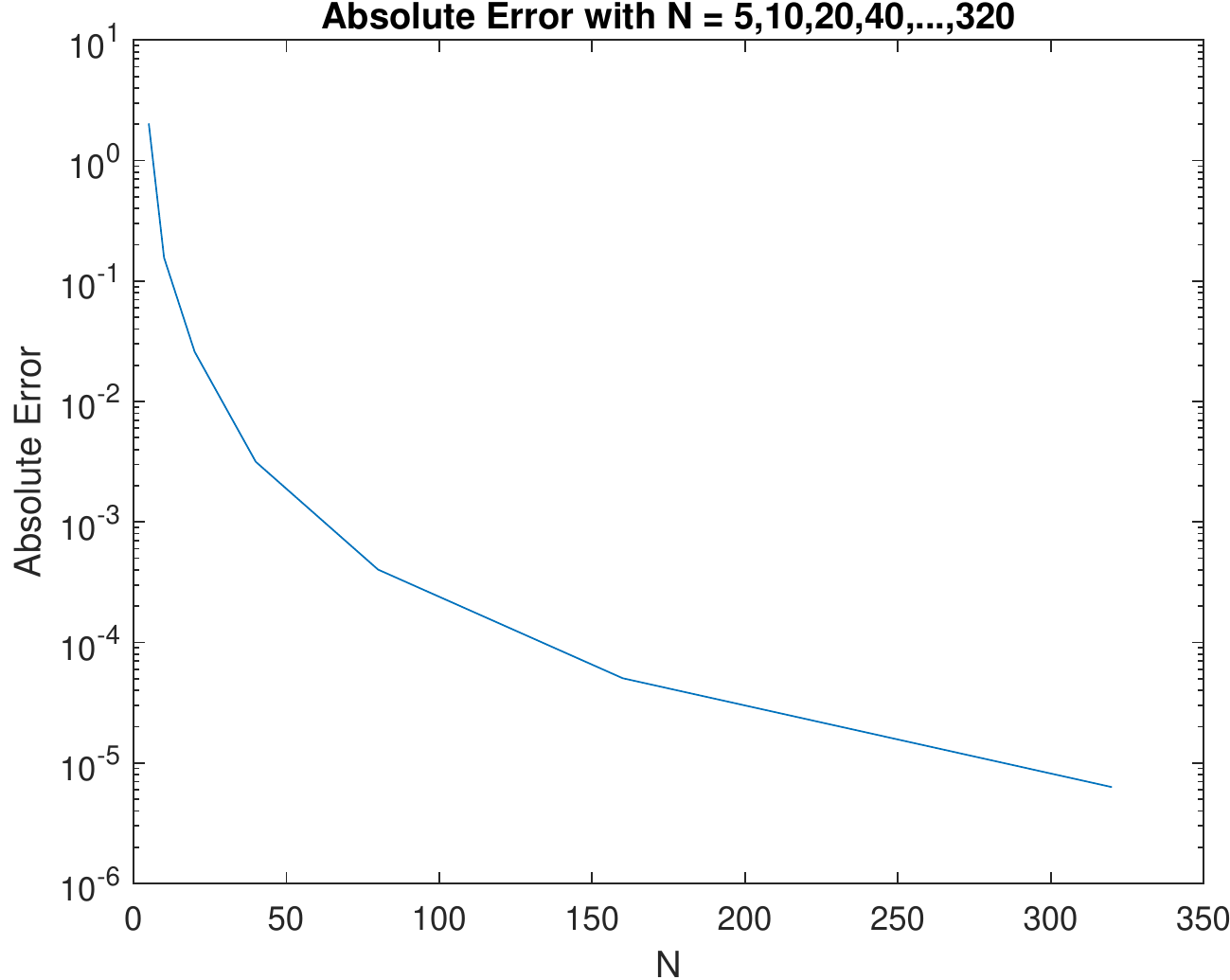}}
     ~ 
    \subfigure[]
        {\includegraphics[height=1.80in]{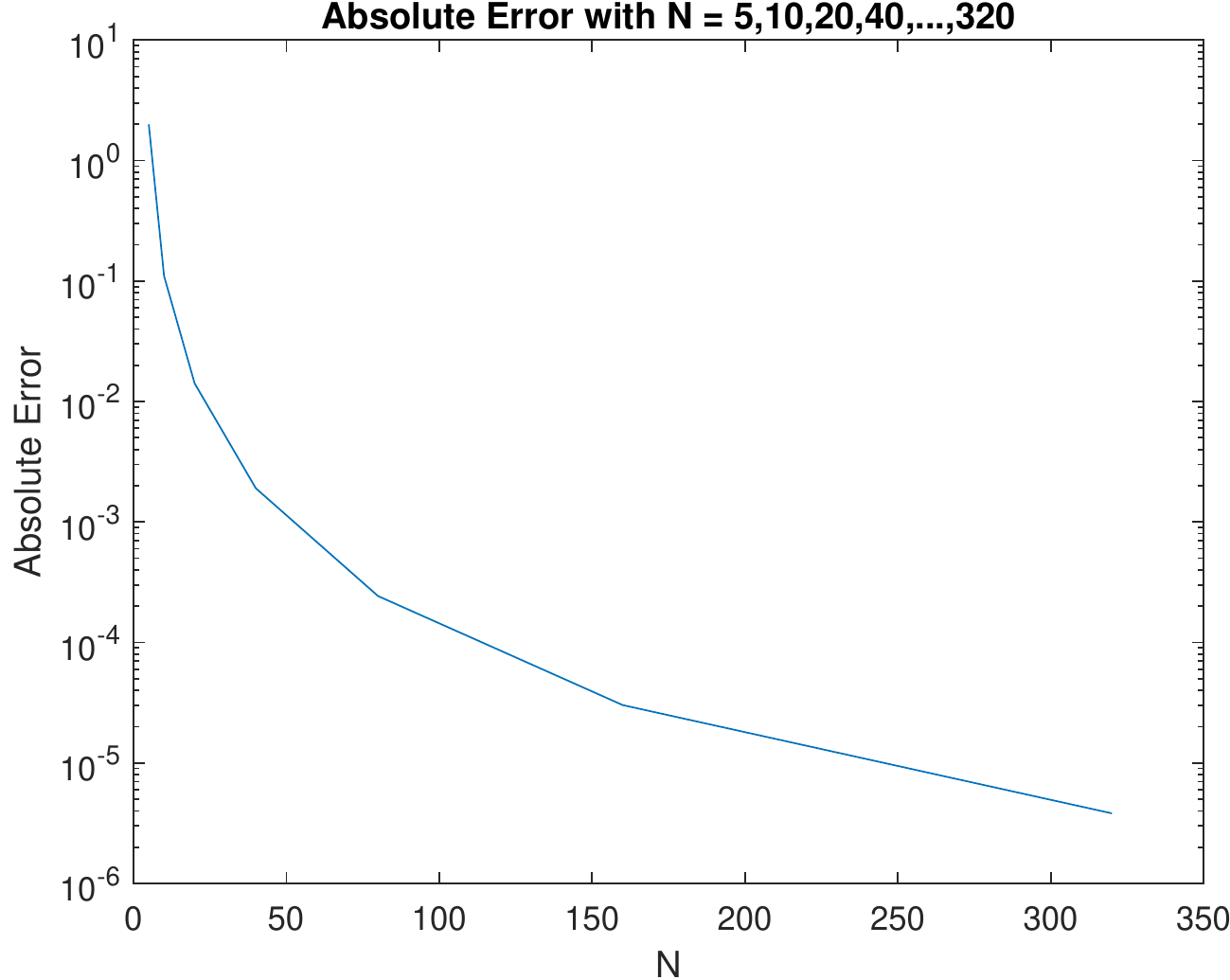}}
    \caption{Convergence of $t_3$ for $f_1$ (left) and $f_2$ (right) and $N=5,10,20,40,80,160,320$.    }
    \label{Fig_Conv_H3}
\end{figure}

\begin{figure}
    \centering
    \graphicspath{ {Images/}}
    \subfigure[]
        {\includegraphics[height=1.80in]{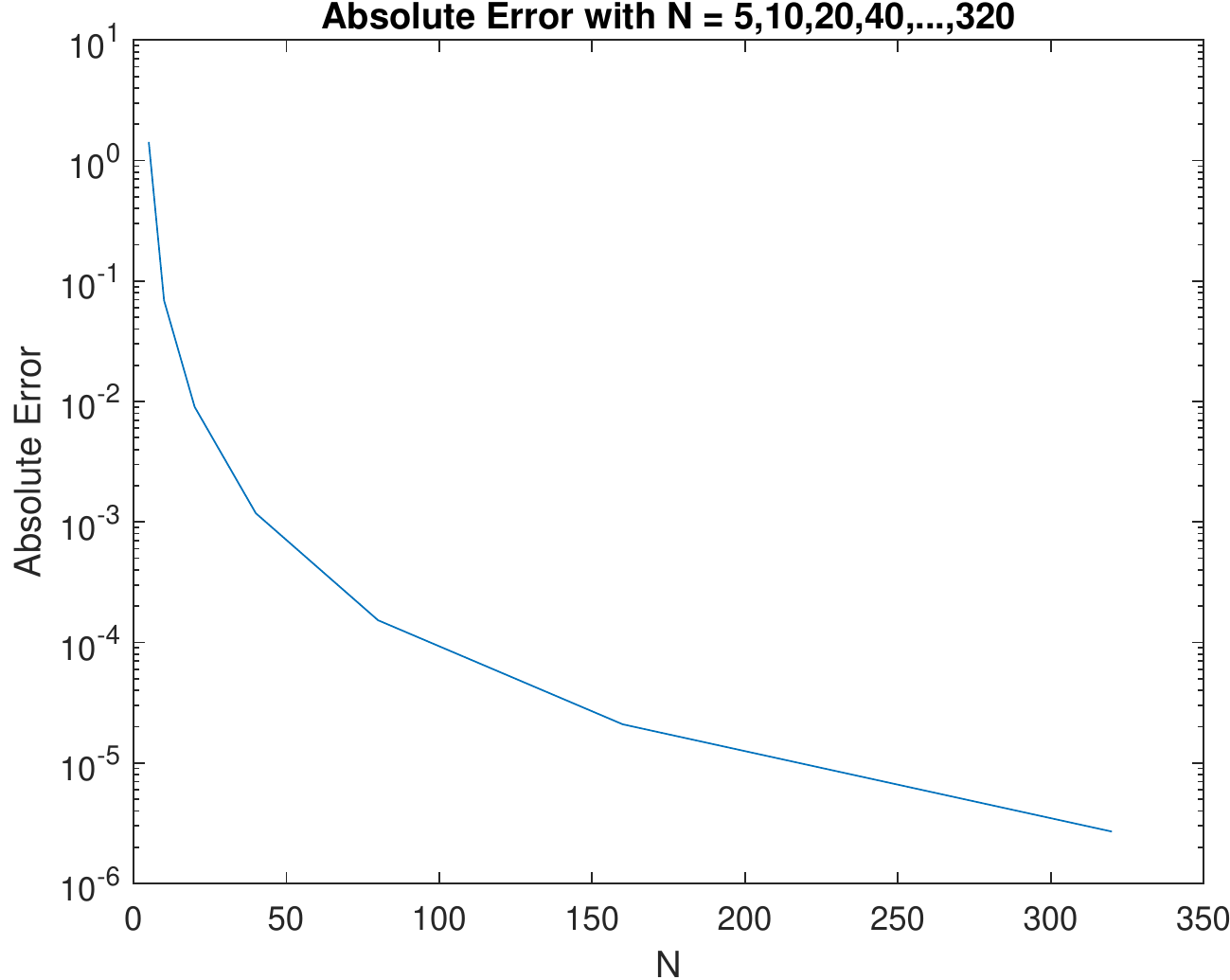}}
     ~ 
    \subfigure[]
        {\includegraphics[height=1.80in]{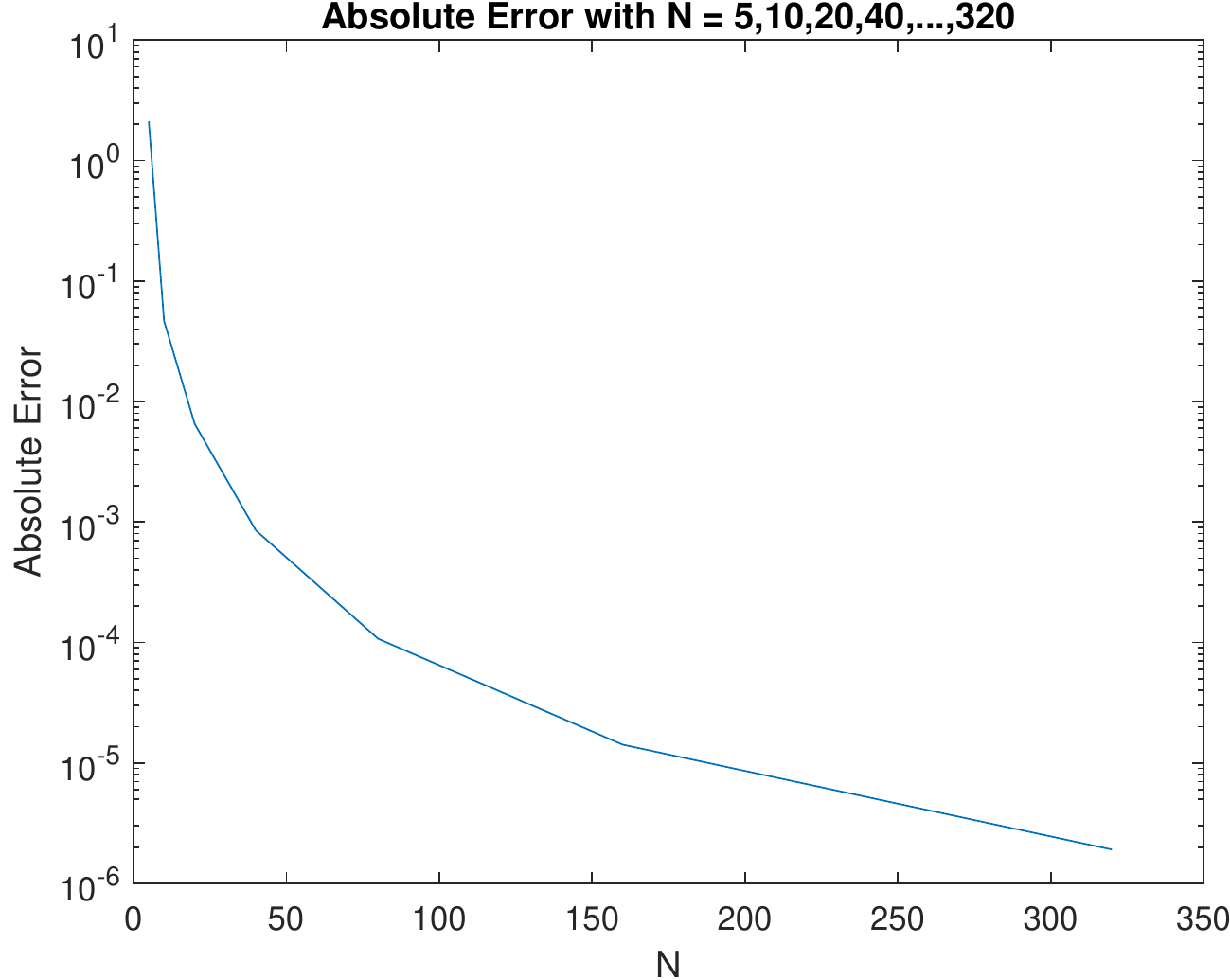}}
    \caption{Convergence of $t_4$ for $f_1$ (left) and $f_2$ (right) and $N=5,10,20,40,80,160,320$.    }
    \label{Fig_Conv_H4}
\end{figure}

For a fixed $m$, we compute, the absolute error between the function and the interpolant with $N$ nodes $\textrm{err}_N = \max_{x\in[0,2\pi)} | f(x) - t_m(x) |$ and we can notice that the interpolant, in the case of conditions up to the second derivative retains the exponential convergence of the classical interpolant \eqref{BaryTrigEqui} at equidistant nodes (see Figures \ref{Fig_Conv_H1} and \ref{Fig_Conv_H2}), on the other hand, when we have conditions on the third or the fourth derivatives, the convergence slows down (see Figures \ref{Fig_Conv_H3} and \ref{Fig_Conv_H4}). We can numerically estimate the convergence of the interpolant with $N$ nodes as $\calO(N^{-3})$. In fact, once we compute $-\log_2(\frac{\textrm{err}_{2n}}{\textrm{err}_{n}})$ with various $n$, we have a numerical estimate of the rate of convergence and it appears, as we can see in Table \ref{Tab_Conv_Hx} where we computed it for $n=5,10,20,40,80,160$, that this ratio is almost $3$.

\begin{table}
    \centering
    \begin{tabular}{l|c||c|c|c|c|c|c}
    & n  &  5 & 10 & 20 & 40 & 80 & 160 \\ 
    \hline
    \multirow{2}{*}{$t_3$} & $f_1$ & 3.70 & 2.59 & 3.03 & 2.97 & 3.00 & 3.00 \\
    & $f_2$ & 4.18 & 2.96 & 2.89 & 2.97 & 3.01 & 2.99 \\
    \hline
    \multirow{2}{*}{$t_4$} & $f_1$ & 4.37 & 2.93 & 2.93 & 2.95 & 2.87 & 2.96 \\
    & $f_2$ & 5.50 & 2.84 & 2.94 & 2.99 & 2.92 & 2.89
    \end{tabular}
    \caption{Computation of the values $-\log_2(\frac{\textrm{err}_{2n}}{\textrm{err}_{n}})$ for $t_3$ and $t_4$ with both test functions}
    \label{Tab_Conv_Hx}
\end{table}

In order to test the interpolant with conformally shifted nodes and explore the convergence to functions with fronts, we tested the interpolant moreover with the function 
$$ f_3(\theta) = \tanh(50\cos(\theta+\pi/3));$$
which has two fronts located in $\theta_1 = \pi/6$ and $\theta_2 = 7\pi/6$. For this reason we compared the performances of the interpolant with equidistant nodes and conformally shifted nodes, via the conformal map introduced in \cite{BE20}, i.e., 
$$ g_{\alpha,\beta} (\theta) = - \iota \log \Bigg( \frac{e^{\iota \theta } + \alpha e^{\iota \beta} }{1 + \alpha e^{\iota (\theta-\beta)}}\Bigg),$$
and using the generalisation for two fronts, where we set the same density parameter for both fronts as $\alpha=0.85$.

\begin{figure}
    \centering
    \graphicspath{ {Images/}}
    \subfigure[]
        {\includegraphics[height=1.75in]{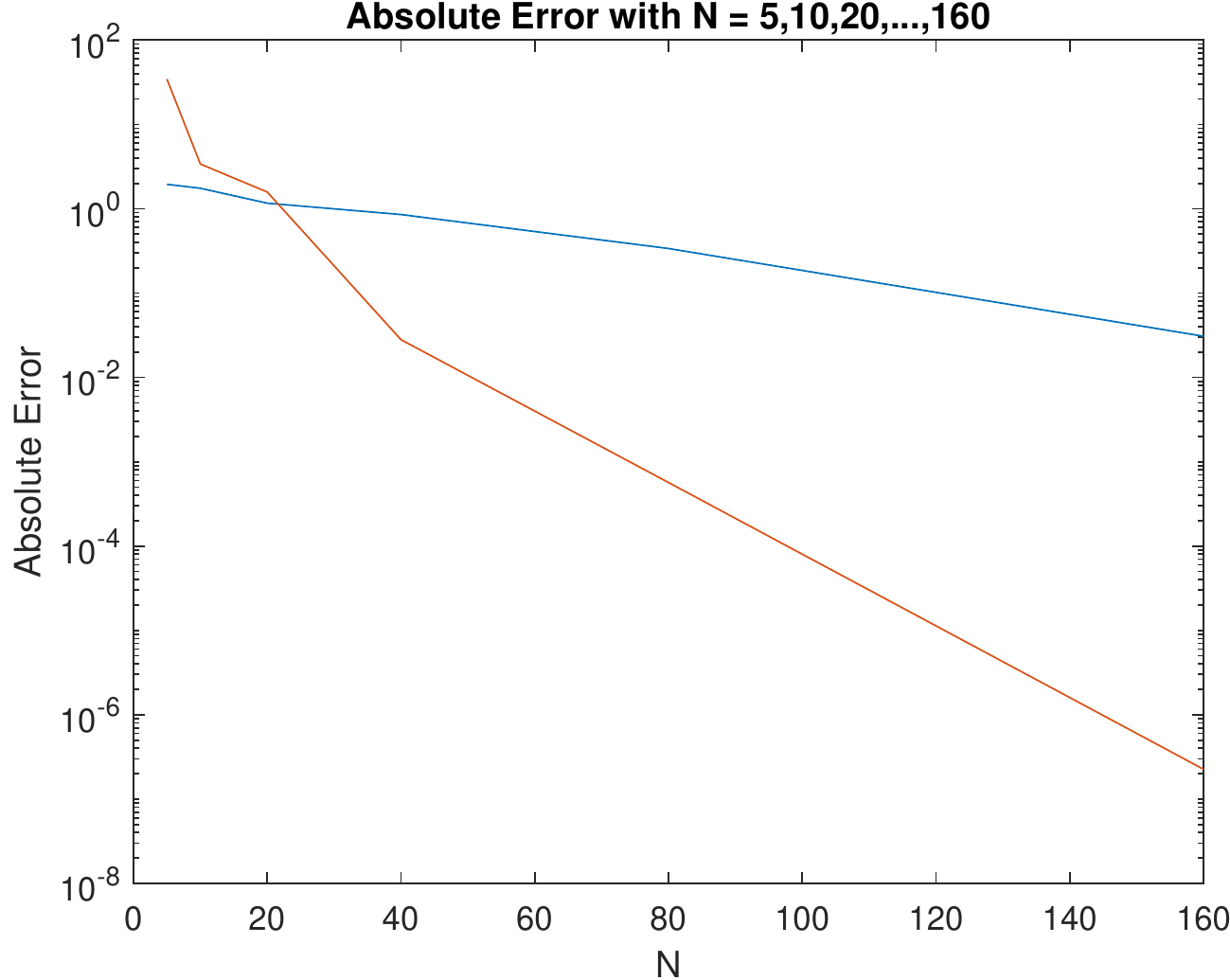}}
     ~ 
    \subfigure[]
        {\includegraphics[height=1.75in]{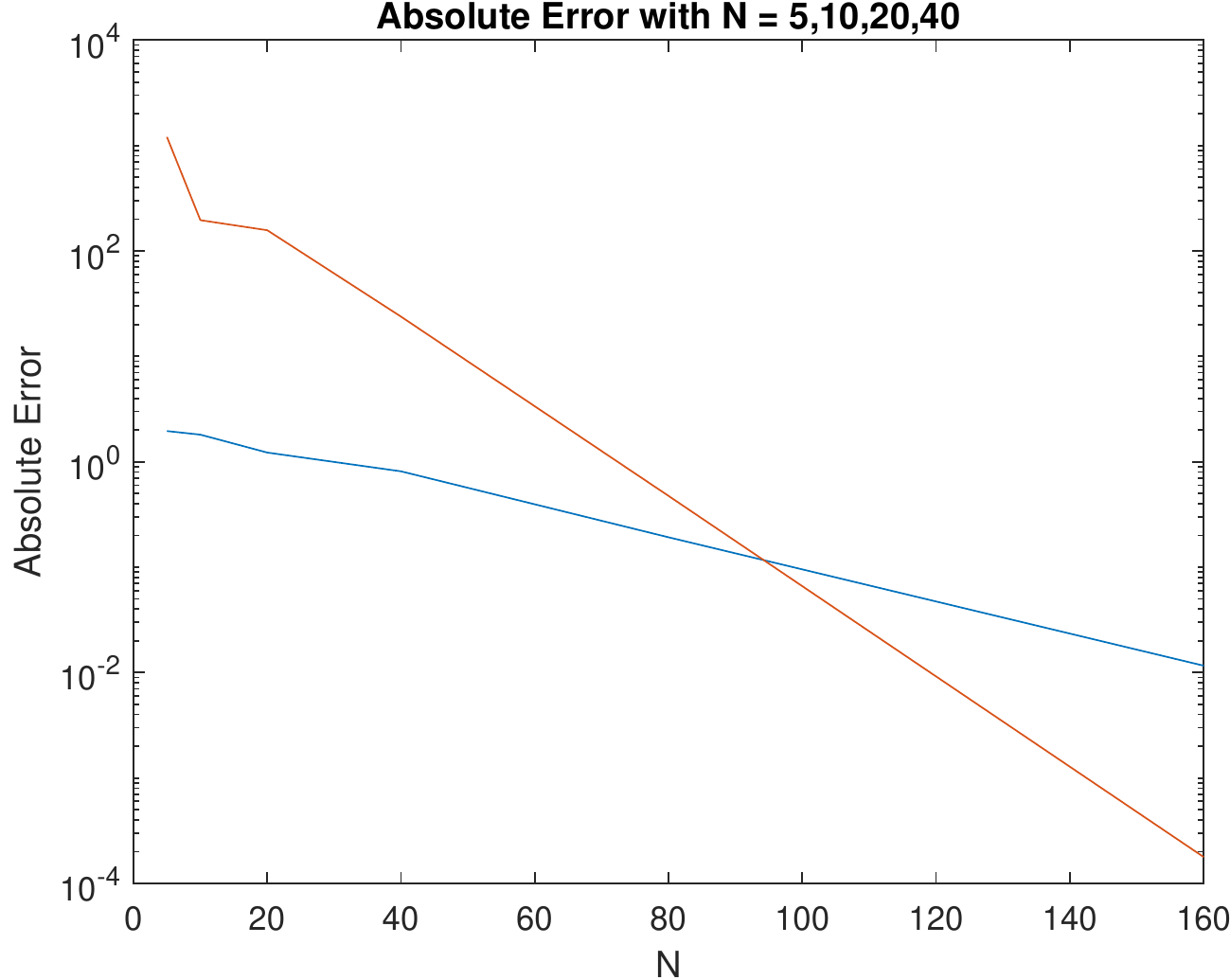}}
     ~ 
    \subfigure[]
        {\includegraphics[height=1.75in]{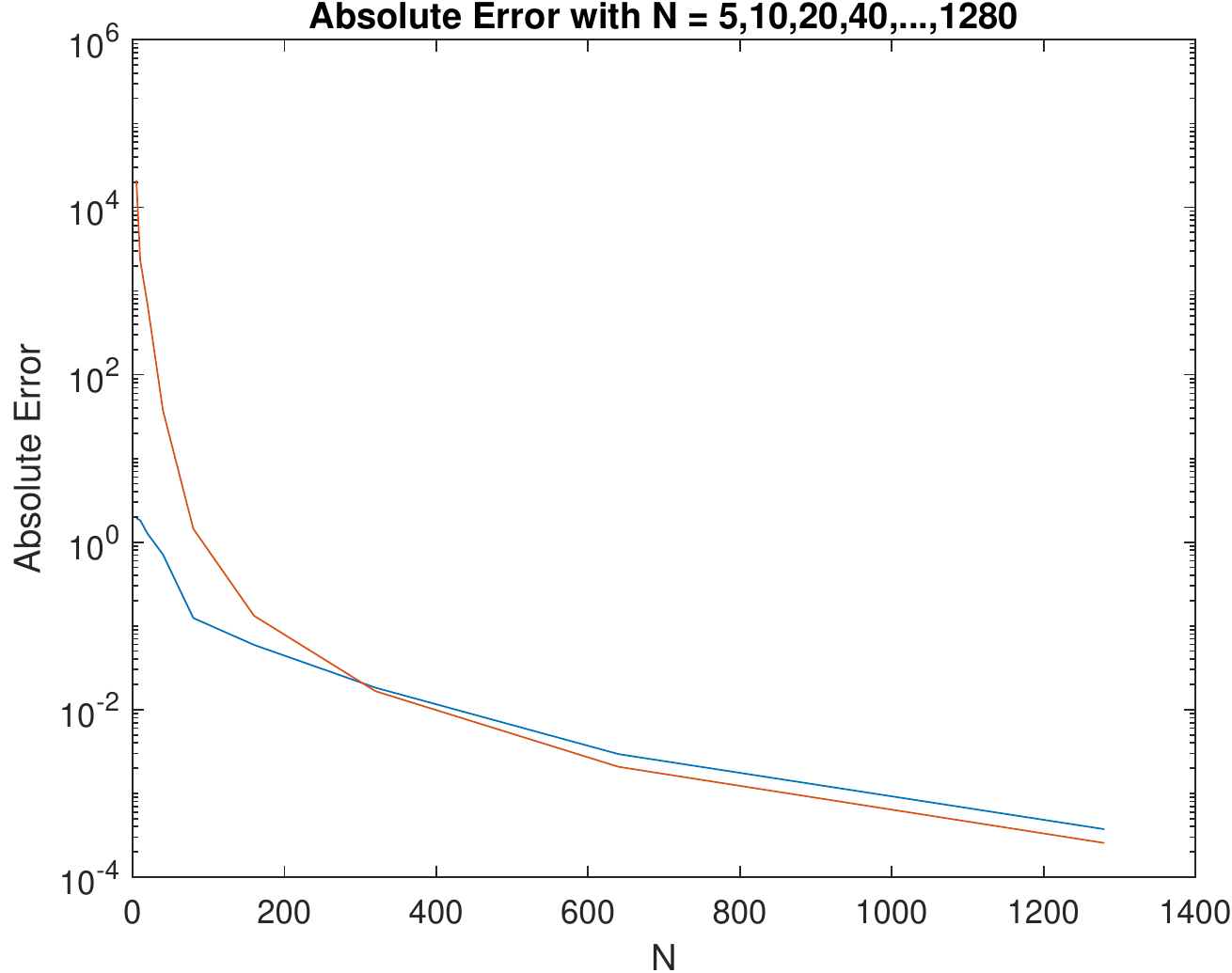}}
     ~ 
    \subfigure[]
        {\includegraphics[height=1.75in]{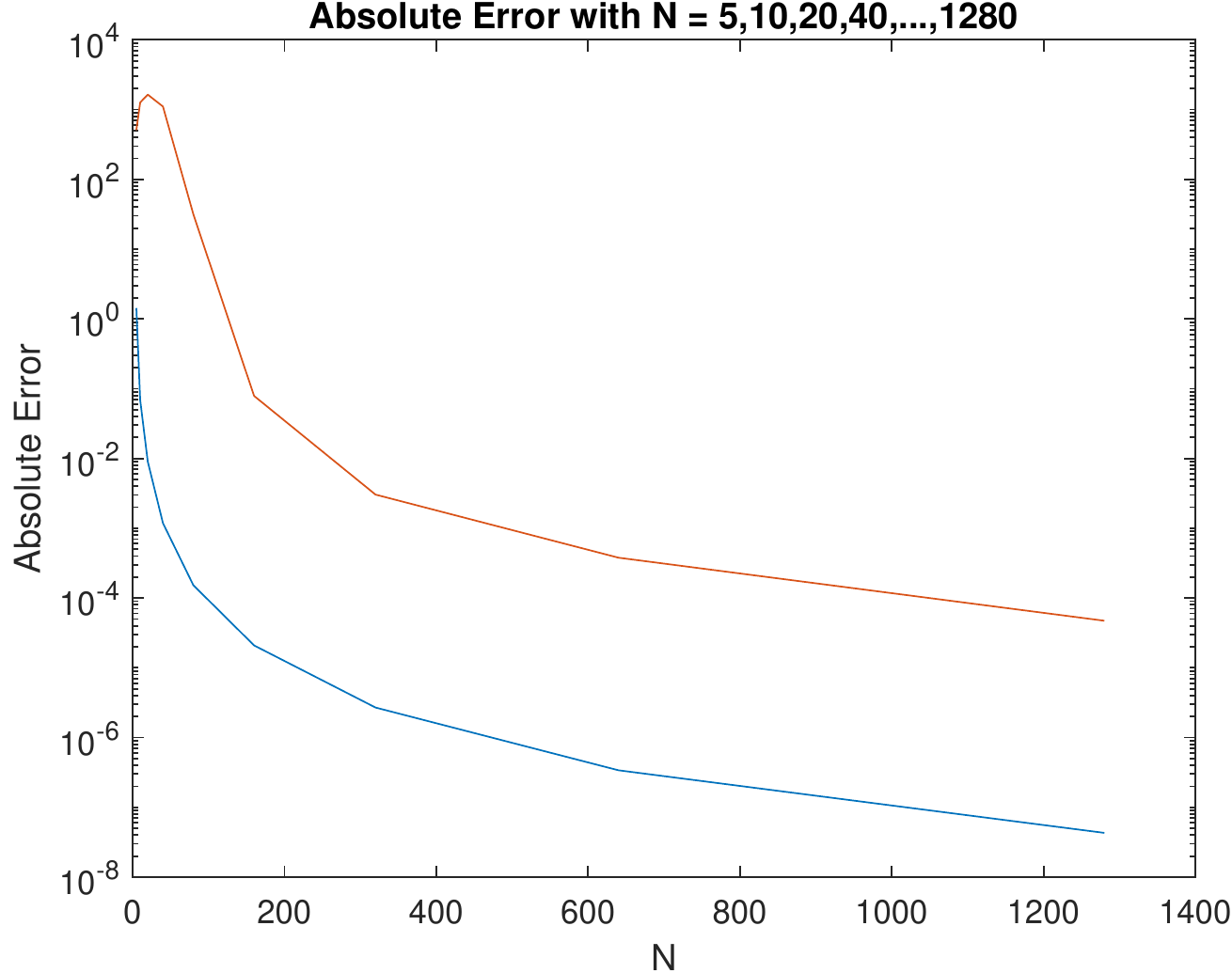}}
    \caption{Error of function $f_3$ with equidistant nodes (blue line) and conformally shifted nodes (red line) by using (a) $t_1$, (b) $t_2$, (c) $t_3$, (d) $t_4$, with $N=5,10,20,\dots,160$ in (a) and (b) and $N=5,10,20,\dots,1280$ in (c) and (d) }
    \label{ConformalShift}
\end{figure}


As we can see in Figures \ref{ConformalShift}, the advantage of clustering the nodes in the locations of the fronts becomes less important as we consider more derivatives. In fact, it speeds up a lot the convergence in the case of one derivative but, then, once we add the conditions on more derivatives, clustering the nodes decreases the speed up until it even slows the convergence when we have conditions up to the fourth derivative. 

We present another numerical experiment which is the reconstruction of a closed simple curve, given some points, as might be a figure or a sketch. The two dimensions are reconstructed separately and each coordinates is consider as a periodic function. Therefore, in order to enjoy the fast exponential convergence, the data are considered as the values of a function at equidistant nodes. Then, we approximate the first derivatives by using the coordinates of the difference quotients of the points.

\begin{figure}
    \centering
    \graphicspath{ {Images/}}
    \subfigure[]
        {\includegraphics[height=1.75in]{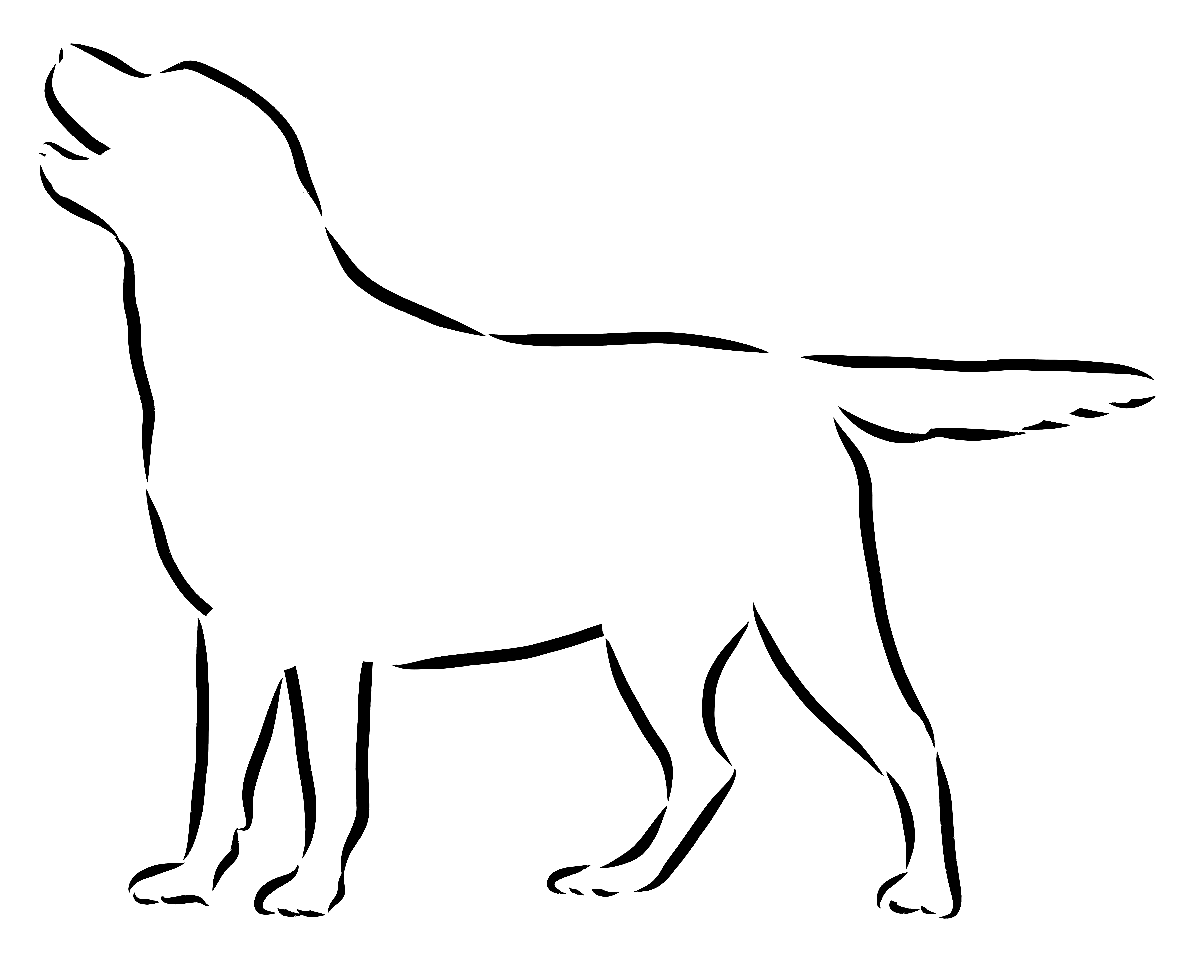}}
     ~ 
    \subfigure[]
        {\includegraphics[height=1.75in]{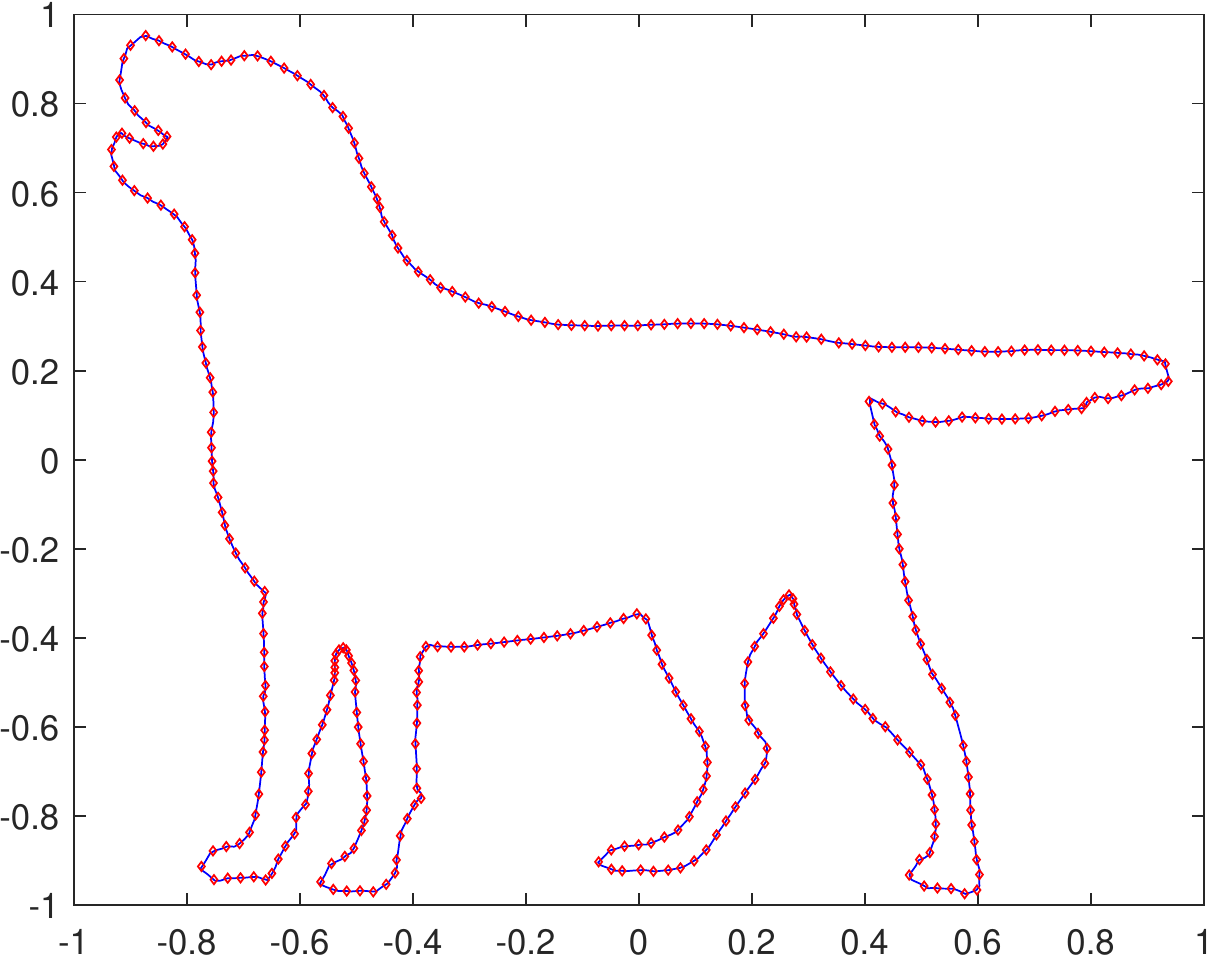}}
     ~ 
    \subfigure[]
        {\includegraphics[height=1.75in]{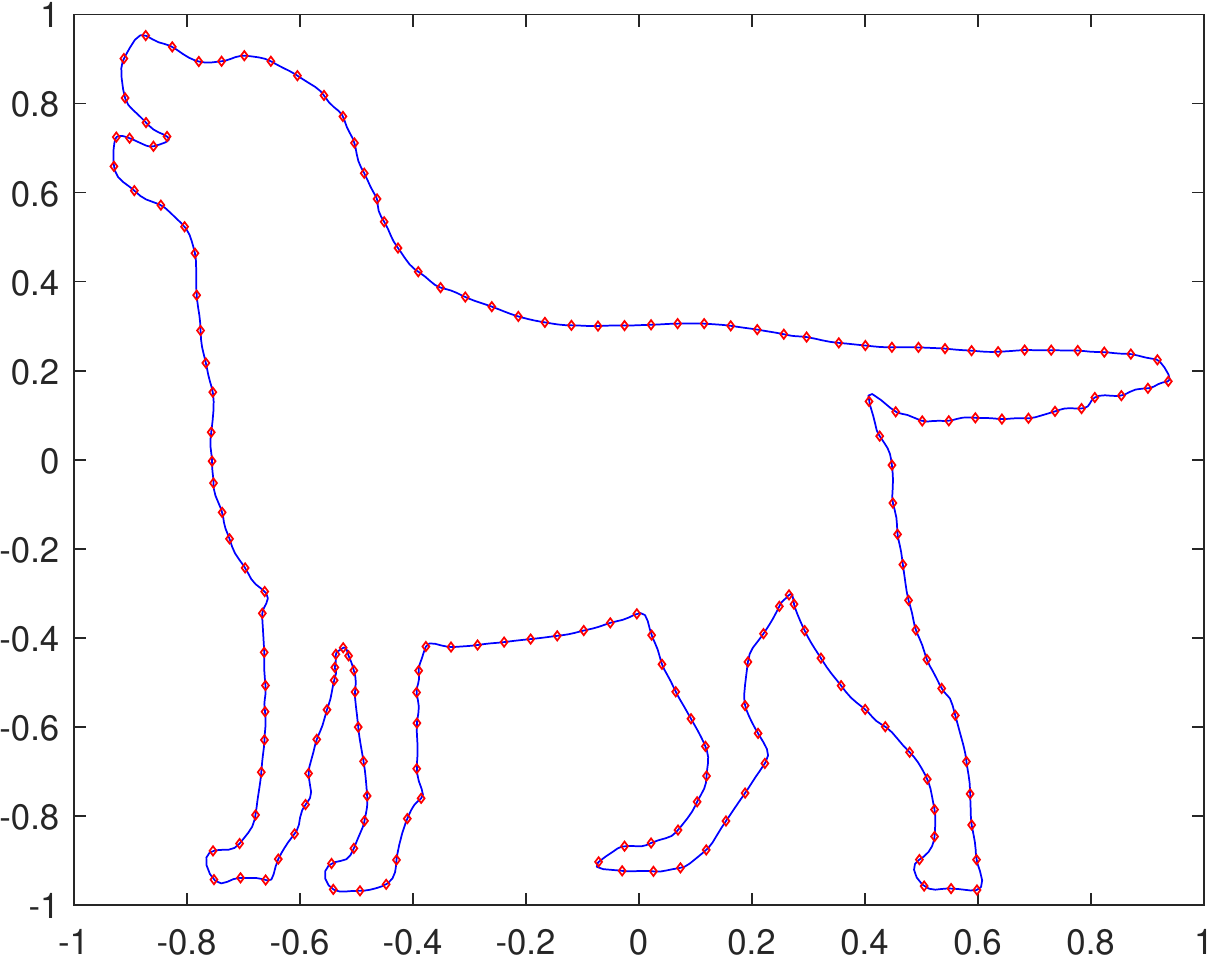}}
     ~ 
    \subfigure[]
        {\includegraphics[height=1.75in]{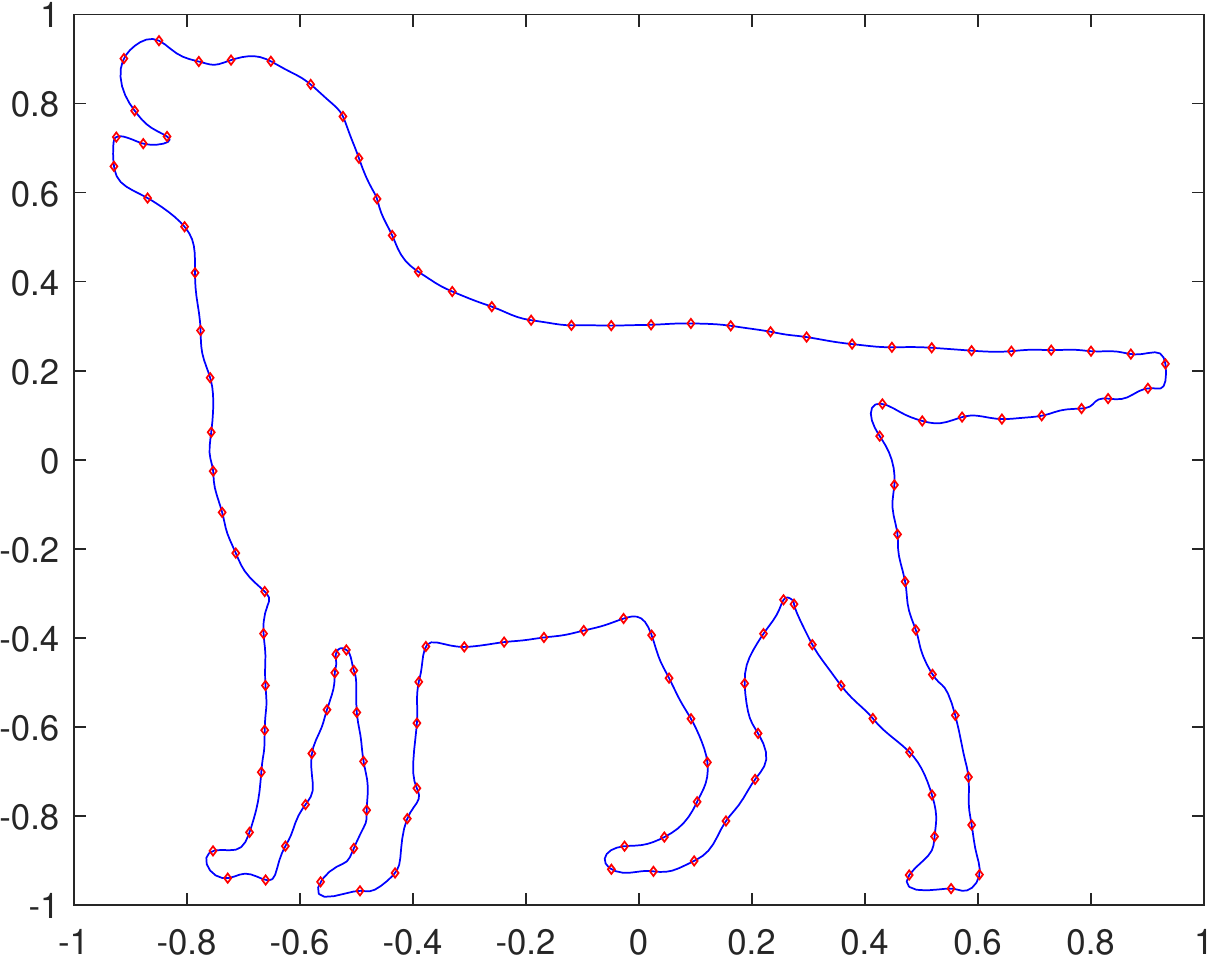}}
     ~ 
    \subfigure[]
        {\includegraphics[height=1.75in]{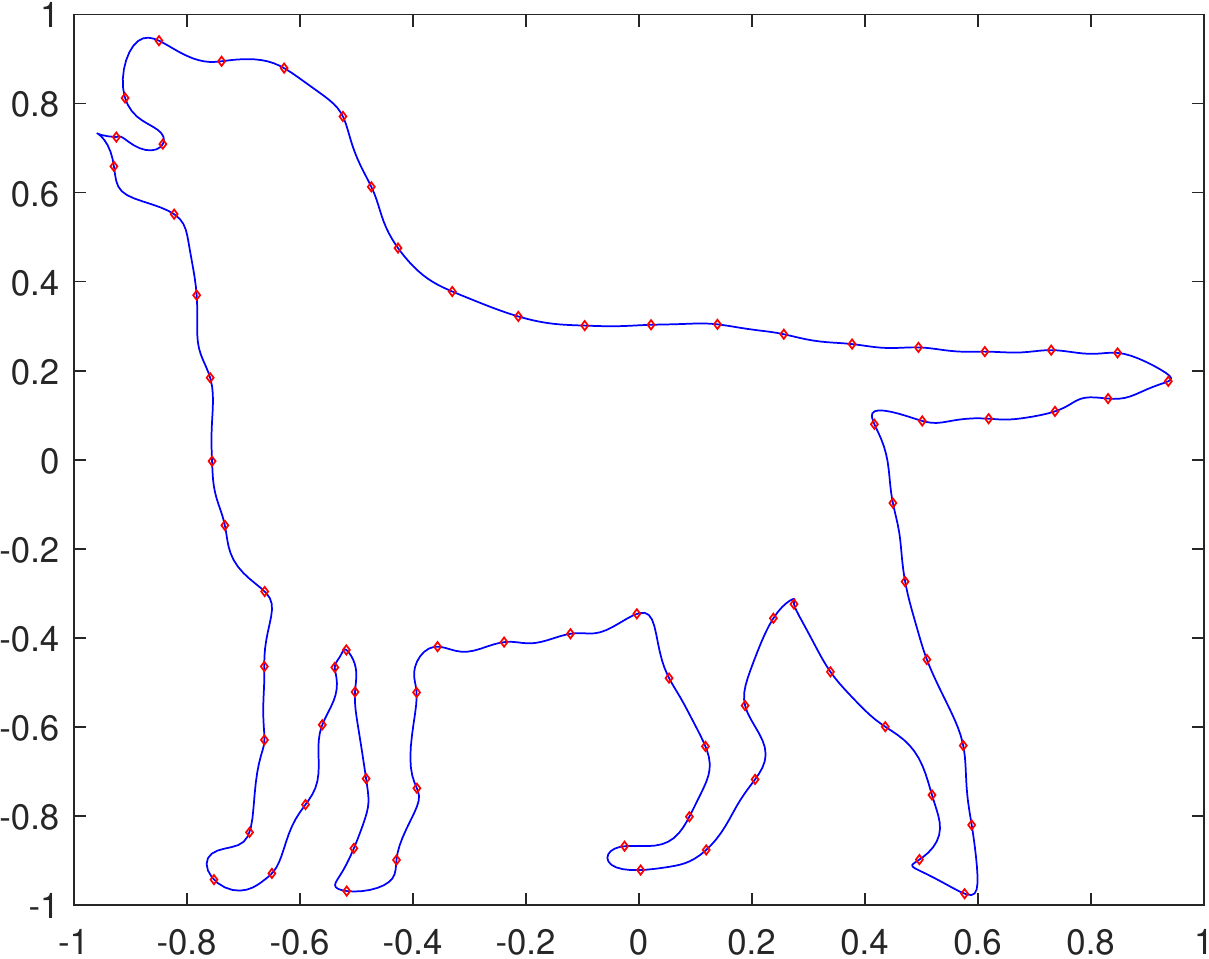}}
    \caption{Reconstruction of the sketch (a), by using 362 nodes (b), 181 nodes (c), 121 nodes (d) and 73 nodes (e).}
    \label{Doggo}
\end{figure}

We did therefore, the reconstruction of the sketch in Figure \ref{Doggo}a and starting by extracting 362 points (which are the red dots in \ref{Doggo}b ), we used the values of the abscissae and ordinates as datae of two periodic function with those values at equidistant nodes and which we intend to approximate by the Hermite interpolant $t_1$, which is the interpolant with conditions on one derivative. Moreover, since the datae of the real tangent vectors are unknown, we approximate them with the difference quotient vector between the extracted points ordered in anticlockwise sense. The result is in Figure \ref{Doggo}b.

In Figure \ref{Doggo}(c)-(e) we used a subset of the initial set of nodes: half points, a third and a fifth, respectively.

As we can see by the results in Figures \ref{Doggo}, the interpolant performs quite well even by using an approximation of the derivative and by decreasing the number of nodes.

In particular, even by using a fifth of the points initially extracted the reconstruction is still resembling quite a lot the sketch.

\section{Conclusion}

The present work introduced an iterative scheme to construct an Hermite interpolant starting from a general interpolation basis and in particular, we focused on Berrut's trigonometric interpolant. In order to construct the values of the derivative of the precedent iterative step, we compute analytically the differentation matrix which allows to obtain these values easily as a matrix product.

Our numerical experiments suggest that the produced interpolant retains the exponential convergence of the classical interpolant at equidistant nodes and conformally shifted nodes, once we interpolate the first two derivatives and slows down to a rate of $\calO(N^{-3})$, where $N$ is the number of nodes used, when we interpolate also the third or the fourth derivatives.

Moreover, the tests suggest that clustering nodes in the locations of fronts via a conformal map speeds up the convergence in case of one derivative but, instead, could be of negative impact when we need to satisfy conditions also for higher derivatives. 

Future work will be aimed into analysing the convergence more deeply, also for other class of nodes, as the one introduced in \cite{BE20}, which includes conformally shifted nodes, and also into analysing the Lebesgue constant of the interpolant.

\vskip 0.1in
{\bf Acknowledgments.} The author thank the referees for their careful reading
of the manuscript and their numerous suggestions, which have improved this work.
This research has been accomplished within Rete ITaliana di Approssimazione (RITA), the thematic group on \lq\lq Teoria dell'Ap\-pros\-si\-ma\-zio\-ne e Applicazioni\rq\rq $\,$(TAA) of the Italian Mathematical Union and partially supported by GNCS-IN$\delta$AM \lq\lq Finanziamenti Giovani Ricercatori 2022\rq\rq $\,$(Young Researchers Fund 2022) for the project \lq\lq An iterative approach for Hermite interpolants\rq\rq,  INdAM - GNCS Project cod. CUP\textunderscore E55F22000270001.


\printbibliography

@article{CH18,
    AUTHOR = {Cirillo, E. and Hormann, K.},
     TITLE = {An iterative approach to barycentric rational {H}ermite
              interpolation},
   JOURNAL = {Numer. Math.},
  FJOURNAL = {Numerische Mathematik},
    VOLUME = {140},
      YEAR = {2018},
    NUMBER = {4},
     PAGES = {939--962},
}

@article{Baltensperger02,
    AUTHOR = {Baltensperger, R.},
     TITLE = {Some results on linear rational trigonometric interpolation},
   JOURNAL = {Comput. Math. Appl.},
  FJOURNAL = {Computers \& Mathematics with Applications. An International
              Journal},
    VOLUME = {43},
      YEAR = {2002},
    NUMBER = {6-7},
     PAGES = {737--746},
}

@article{Berrut88,
    AUTHOR = {Berrut, J.-P.},
     TITLE = {Rational functions for guaranteed and experimentally
              well-conditioned global interpolation},
   JOURNAL = {Comput. Math. Appl.},
  FJOURNAL = {Computers \& Mathematics with Applications. An International
              Journal},
    VOLUME = {15},
      YEAR = {1988},
    NUMBER = {1},
     PAGES = {1--16},
}

@article {BE21,
    AUTHOR = {Berrut, J.-P. and Elefante, G.},
     TITLE = {Bounding the {L}ebesgue constant for a barycentric rational
              trigonometric interpolant at periodic well-spaced nodes},
   JOURNAL = {J. Comput. Appl. Math.},
  FJOURNAL = {Journal of Computational and Applied Mathematics},
    VOLUME = {398},
      YEAR = {2021},
     PAGES = {113664, 11},
}

@article {Welfert97,
    AUTHOR = {Welfert, B.~D.},
     TITLE = {Generation of pseudospectral differentiation matrices. {I}},
   JOURNAL = {SIAM J. Numer. Anal.},
  FJOURNAL = {SIAM Journal on Numerical Analysis},
    VOLUME = {34},
      YEAR = {1997},
    NUMBER = {4},
     PAGES = {1640--1657},
}

@article {BE20,
    AUTHOR = {Berrut, J.-P. and Elefante, G.},
     TITLE = {A periodic map for linear barycentric rational trigonometric
              interpolation},
   JOURNAL = {Appl. Math. Comput.},
  FJOURNAL = {Applied Mathematics and Computation},
    VOLUME = {371},
      YEAR = {2020},
     PAGES = {124924, 8},
}

@article {CH19,
    AUTHOR = {Cirillo, E. and Hormann, K.},
     TITLE = {On the {L}ebesgue constant of barycentric rational {H}ermite
              interpolants at equidistant nodes},
   JOURNAL = {J. Comput. Appl. Math.},
  FJOURNAL = {Journal of Computational and Applied Mathematics},
    VOLUME = {349},
      YEAR = {2019},
     PAGES = {292--301},
}

@article {CH20,
    AUTHOR = {Cirillo, E. and Hormann, K. and Sidon, J.},
     TITLE = {Convergence rates of a {H}ermite generalization of
              {F}loater-{H}ormann interpolants},
   JOURNAL = {J. Comput. Appl. Math.},
  FJOURNAL = {Journal of Computational and Applied Mathematics},
    VOLUME = {371},
      YEAR = {2020},
     PAGES = {112624, 9},
}

@article {BD19,
    AUTHOR = {Bandiziol, C. and De Marchi, S.},
     TITLE = {On the {L}ebesgue constant of the trigonometric
              {F}loater-{H}ormann rational interpolant at equally spaced
              nodes},
   JOURNAL = {Dolomites Res. Notes Approx.},
  FJOURNAL = {Dolomites Research Notes on Approximation},
    VOLUME = {12},
      YEAR = {2019},
     PAGES = {51--67},
}

@Article{Oumellal21,
AUTHOR = {Oumellal, F. and Lamnii, A.},
TITLE = {Curve and Surface Construction Using Hermite Trigonometric Interpolant},
JOURNAL = {Mathematical and Computational Applications},
VOLUME = {26},
YEAR = {2021},
NUMBER = {1},
ARTICLE-NUMBER = {11},
}

@article{Jing20,
    AUTHOR = {Jing, K. and Liu, Y. and Kang, N. and Zhu, G.},
     TITLE = {A convergent family of linear {H}ermite barycentric rational
              interpolants},
   JOURNAL = {J. Math. Res. Appl.},
  FJOURNAL = {Journal of Mathematical Research with Applications},
    VOLUME = {40},
      YEAR = {2020},
    NUMBER = {6},
     PAGES = {628--646},
}

@article{Krasny19,
    AUTHOR = {Krasny, R. and Wang, L.},
     TITLE = {A treecode based on barycentric {H}ermite interpolation for
              electrostatic particle interactions},
   JOURNAL = {Comput. Math. Biophys.},
  FJOURNAL = {Computational and Mathematical Biophysics},
    VOLUME = {7},
      YEAR = {2019},
     PAGES = {73--84},
}

@article{Berrut22,
    AUTHOR = {Berrut, J.-P. and Elefante, G.},
     TITLE = {A linear barycentric rational interpolant on starlike domains},
   JOURNAL = {Electron. Trans. Numer. Anal.},
  FJOURNAL = {Electronic Transactions on Numerical Analysis},
    VOLUME = {55},
      YEAR = {2022},
     PAGES = {726--743},
}

@article{Schneider91,
    AUTHOR = {Schneider, C. and Werner, W.},
     TITLE = {Hermite interpolation: the barycentric approach},
   JOURNAL = {Computing},
  FJOURNAL = {Computing. Archives for Scientific Computing},
    VOLUME = {46},
      YEAR = {1991},
    NUMBER = {1},
     PAGES = {35--51},
}

@article {DellAccio20,
    AUTHOR = {Dell'Accio, F. and Di Tommaso, F. and Nouisser,
              O. and Siar, N.},
     TITLE = {Rational {H}ermite interpolation on six-tuples and scattered data},
   JOURNAL = {Appl. Math. Comput.},
  FJOURNAL = {Applied Mathematics and Computation},
    VOLUME = {386},
      YEAR = {2020},
     PAGES = {125452, 11},
}

@article {Daus20,
    AUTHOR = {D\u{a}u\c{s}, L. and Jianu, M.},
     TITLE = {Full {H}ermite interpolation of the reliability of a hammock
              network},
   JOURNAL = {Appl. Anal. Discrete Math.},
  FJOURNAL = {Applicable Analysis and Discrete Mathematics},
    VOLUME = {14},
      YEAR = {2020},
    NUMBER = {1},
     PAGES = {198--220},
}

@article{Henrici79,
    AUTHOR = {Henrici, P.},
     TITLE = {Barycentric formulas for interpolating trigonometric
              polynomials and their conjugates},
   JOURNAL = {Numer. Math.},
  FJOURNAL = {Numerische Mathematik},
    VOLUME = {33},
      YEAR = {1979},
    NUMBER = {2},
     PAGES = {225--234},
}

@article {Qi15,
    AUTHOR = {Qi, F.},
     TITLE = {Derivatives of tangent function and tangent numbers},
   JOURNAL = {Appl. Math. Comput.},
  FJOURNAL = {Applied Mathematics and Computation},
    VOLUME = {268},
      YEAR = {2015},
     PAGES = {844--858},
}

@article {Baddoo21,
    AUTHOR = {Baddoo, P. J.},
     TITLE = {The {AAA}trig algorithm for rational approximation of periodic
              functions},
   JOURNAL = {SIAM J. Sci. Comput.},
  FJOURNAL = {SIAM Journal on Scientific Computing},
    VOLUME = {43},
      YEAR = {2021},
    NUMBER = {5},
     PAGES = {A3372--A3392},
}

@article {Nakatsukasa18,
    AUTHOR = {Nakatsukasa, Y. and S\`ete, O. and Trefethen, L.~N.},
     TITLE = {The {AAA} algorithm for rational approximation},
   JOURNAL = {SIAM J. Sci. Comput.},
  FJOURNAL = {SIAM Journal on Scientific Computing},
    VOLUME = {40},
      YEAR = {2018},
    NUMBER = {3},
     PAGES = {A1494--A1522},
}

@article {Baddoo19,
    AUTHOR = {Baddoo, P. J. and Crowdy, D. G.},
     TITLE = {Periodic {S}chwarz-{C}hristoffel mappings with multiple
              boundaries per period},
   JOURNAL = {Proc. A.},
  FJOURNAL = {Proceedings A},
    VOLUME = {475},
      YEAR = {2019},
    NUMBER = {2228},
     PAGES = {20190225, 20},
}

\end{document}